\documentclass[smallextended]{svjour3} 
\smartqed

\usepackage{amssymb}
\usepackage{amsmath}
\usepackage{cases}
\usepackage{enumerate}
\usepackage{graphicx}
\usepackage[sort&compress,numbers]{natbib}
\usepackage{indentfirst}
\usepackage{doi}
\usepackage[usenames]{color}
\usepackage{wrapfig}
\usepackage{txfonts}
\usepackage{stmaryrd}
\usepackage{array} 
\usepackage[caption=false]{subfig} 
\SetSymbolFont{stmry}{bold}{U}{stmry}{m}{n}

\usepackage{mathrsfs}
\usepackage{pgfplots}
\usepackage{caption}
\usepackage{enumitem}
\pgfplotsset{compat=1.7}

\usepackage{latexsym}
\usepackage[numbers]{natbib}
\usepackage{lmodern}

\newcommand*{\vv}[1]{\vec{\mkern0mu#1}}

\newcommand{\norm}[1]{\Vert#1\Vert}

\newcommand{\bR}{{\mathbb R}}

\newcommand{\bZ}{\mathbb{Z}}

\newcommand{\Gt}{{\Gamma(t)}}

\newcommand{\rd}{\;{\rm d}}

\newcommand{\dd}[1]{\frac{\rm d}{{\rm d}#1}}
\newcommand{\ddt}{\dd{t}}
\newcommand{\nn}{\nonumber}
\newcommand{\ttau}{\Delta t}

\journalname{J. Sci. Comput.}

\begin{document}
\title{Stable fully discrete finite element methods  with BGN tangential motion  for Willmore flow of planar curves}
\titlerunning{Stable FEM for Willmore flow of planar curves}  

\author{\makebox[\textwidth][l]{Harald Garcke \and  Robert N\"urnberg  \and   Quan~Zhao}}
\authorrunning{H. Garcke, R. N\"urnberg and Q. Zhao}

\institute{Harald Garcke \at
              Fakult{\"a}t f{\"u}r Mathematik,
              Universit{\"a}t Regensburg, 
              93040 Regensburg, Germany.\\
              \email{harald.garcke@ur.de}          
           \and
           Robert N\"urnberg \at
            Dipartimento di Mathematica, 
            Universit\`a di Trento,
            38123 Trento, Italy.\\
            \email{robert.nurnberg@unitn.it}
            \and 
            Quan Zhao \at 
            School of Mathematical Sciences, University of Science and Technology of China, 230026 Hefei, Anhui, China.\\
              \email{quanzhao@ustc.edu.cn} 
}

\date{Received: date / Accepted: date}

\maketitle

\begin{abstract}
We propose and analyze stable finite element approximations for Willmore flow of planar curves. The presented schemes are based on a novel weak formulation which combines an evolution equation for curvature with
the curvature formulation originally proposed by Barrett, Garcke and N\"urnberg (BGN) in \cite{BGN07}. 
Under discretization in space with piecewise linear elements this leads to a stable continuous-in-time semidiscrete scheme, which retains the equidistribution property from the BGN methods. Furthermore, two fully discrete schemes can be shown to satisfy unconditional energy stability estimates. Numerical examples are presented to showcase the good properties of the introduced schemes, including an asymptotic equidistribution of vertices.

\end{abstract}

\renewcommand{\thefootnote}{\arabic{footnote}}
\setlength{\parindent}{2em}

\setcounter{equation}{0}

\section{Introduction} \label{sec:intro}
\setlength\parindent{10pt}

Geometric functionals that involve curvature quantities play an important role in geometric analysis and applied mathematics, with broad applications in various fields. One of the prominent examples is the Willmore energy, which is defined as the integral of the square of the mean curvature over a surface \cite{Willmore93}. 
In this paper we focus on planar curves, in which case the Willmore energy
reduces to the integral of the square of the curvature, and is often also
called elastic energy. Applications involving the Willmore energy have been considered in the fields of biomembranes \cite{Canham1970minimum, Helfrich73elastic,Seifert97}, computer graphics \cite{Welch94free,Desbrun99implicit}, materials science \cite{ChenLSWW18,Grinspun08}, 
computational geometry \cite{BretinLO11,GruberA20} and surface
restoration \cite{ClarenzDDRR04,LeeC23}.

Willmore flow is the $L^2$-gradient flow of the Willmore energy, and thus can be used to obtain energy minimizers. The governing equations for the flow are fourth-order and highly nonlinear, thus posing significant challenges for the mathematical and numerical analysis. A crucial aspect for the latter is the development of energy stable discretizations, that allow for the efficient integration of the evolution equation using possibly large time steps. In this paper we will restrict our attention to Willmore flow of curves in the plane, which is often also called elastic flow of (planar) curves. It is our aim to introduce a fully discrete parametric finite element approximation that is unconditionally energy stable, and that in addition has nice mesh properties. In what follows, we are going to review the literature on existing parametric finite element methods for Willmore flow of planar curves. For references to other numerical approaches, as well as to theoretical results for Willmore flow, we refer to the review article \cite{DeckelnickDE05}, and to the introduction of the recent work \cite{BaoL25}.

One of the earliest contributions on the parametric finite element
approximation of Willmore flow for curves was given in \cite{DziukKS02}.
Later a stable semidiscrete continuous-in-time scheme was proposed in 
\cite{DeckelnickD09}, and an error estimate was shown for it.
For a fully discrete variant, on assuming a strict CFL condition, 
error estimates and stability was proved in \cite{Bondarava15}.
All these approaches discretize a purely normal flow, which can lead to very 
nonuniform meshes and coalescence of vertices in practice.
Purely normal discretizations for Willmore flow of surfaces have been
considered in \cite{ClarenzDDRR04,Rusu05,Dziuk08, Duan2021high}.  For example, in \cite{Duan2021high}  a
fully discrete energy-stable scheme was proposed based on the collocation method of Dziuk's formulation \cite{Dziuk08}. For completeness, let us
mention that minimizing movement schemes for Willmore flow of curves and
surfaces, which by definition are unconditionally stable, have been considered
in \cite{OlischlagerR09,BalzaniR12}.
In addition, an unconditionally stable finite element method for
the elastic flow of inextensible curves was proposed and studied in
\cite{Bartels13a}.

Let us now return our focus to methods that discretize a weak formulation of a
system of PDEs describing the Willmore flow of curves.
The BGN schemes from the series of papers
\cite{BGN07,BGN08willmore,curves3d} enjoy nearly
uniform distributions of mesh points in practice. A crucial aspect of these
schemes in the planar case is the spatial discretization of the identity
\begin{equation} \label{eq:bgn}
\varkappa\, \vec\nu = \vec x_{ss},
\end{equation}
where $\varkappa$ and $\vec\nu$ denote the curvature and unit normal of the curve 
$\Gamma = \vec x(\mathbb{I})$. This induces what is now often called a BGN
tangential motion, leading to equidistribution in the semidiscrete case, and to
asymptotically equidistributed fully discrete approximations. We refer to the
recent review article \cite{Barrett20} for more details. We stress that for the
original BGN schemes \cite{BGN07,BGN08willmore,curves3d} no stability estimates are available. However, upon
adapting the ideas from \cite{DeckelnickD09}, a stable 
semidiscrete continuous-in-time BGN scheme was introduced in \cite{pwf}.
For over a decade no further progress on stable approximations of Willmore flow
was made, until very recently Bao and Li introduced a novel scheme in
\cite{BaoL25}. By considering a system of PDEs that includes an evolution equation for 
curvature, they are able to prove energy stability for a fully discrete 
parametric finite element approximation. We notice that incorporating evolution
equations for geometric quantities into the system of PDEs to discretize
is an idea previously exploited by Kov\'acs, Li and Lubich
in \cite{KovacsLL19,KovacsLL21} for mean curvature flow and Willmore flow of
surfaces in $\bR^3$, respectively.
The novel scheme from \cite{BaoL25} exhibits a tangential motion that is not of BGN
type, and as a consequence for complex curve evolutions their
numerical scheme may require mesh redistributions. 

It is our aim to adapt the ideas from \cite{BaoL25} such that the
obtained fully discrete schemes are energy stable, and such that even complicated curve
evolutions can be computed without mesh redistribution. We will accomplish this by
successfully marrying an evolution equation for curvature to the identity
\eqref{eq:bgn}. In fact this can be achieved by utilizing two independent
curvature discretizations, rather than a single one. While one discrete curvature
will approximate the gradient flow structure of the flow, the second
discrete curvature acts like a Lagrange multiplier for the side constraint
\eqref{eq:bgn}.
Finally, let us mention that for all the above schemes that do not 
discretize a purely normal flow, 
a convergence analysis is still missing. The first error analysis
for a parametric finite element approximation of Willmore flow for curves
with nice tangential motion was recently presented in \cite{cd}. There the
authors formulate a strictly parabolic system of PDEs that includes a
well-defined tangential motion, which in practice leads to nearly
equidistributed curves. However, a stability result for the new scheme from
\cite{cd} is not available.

The rest of the paper is organized as follows. In Section~\ref{sec:mathf}, we first derive a new system of geometric PDEs for Willmore flow of planar curves, and then propose a novel weak formulation that allows for BGN tangential movement under discretization. Next in Section~\ref{sec:semi}, we consider a semidiscrete approximation of the weak formulation and prove that the approximation satisfies a stability bound and an equidistribution property. Subsequently in Section~\ref{sec:fem}, we employ temporal discretization to obtain two fully discrete finite element schemes and prove their unconditional energy stability. 
We further report on numerical experiments in Section~\ref{sec:num}, which 
demonstrate the convergence, stability and good mesh properties of our 
introduced schemes. Finally, we draw some conclusions in Section~\ref{sec:con}.

\section{Mathematical formulations}\label{sec:mathf}

\subsection{Problem setting}\label{sec:ps}
Let $\Gamma(t)$ be a closed curve in $\bR^2$ with a smooth parameterization given by 
\begin{equation*}
\vec x(\rho,t): \mathbb{I}\times[0,T]\mapsto\mathbb{R}^2,\nn
\end{equation*}
where $\mathbb{I}=\bR/\bZ$ is the reference domain. On assuming that $|\vec x_\rho|>0$, we introduce the arc length of the curve $\Gamma(t)$, i.e., $\partial_s = |\vec x_\rho|^{-1}\partial_\rho$. We also introduce the unit tangent and normal to the curve $\Gamma(t)$ via 
\begin{equation}\label{eqn:taun}
\vec\tau(\rho,t) = \vec x_s = |\vec x_\rho|^{-1}\vec x_\rho, \qquad \vec\nu(\rho,t) = -(\vec\tau)^\perp=-(\vec x_s)^\perp,
\end{equation}
where $(\cdot)^\perp$ denotes a clockwise rotation by $\frac{\pi}{2}$. The curvature of the curve is defined as 
\begin{equation}\label{eq:cur}
\varkappa(\rho,t) = \vec x_{ss}\cdot\vec\nu\qquad\forall\rho\in\mathbb{I}.
\end{equation}
We further introduce the velocity induced by this parameterization and the normal velocity of the curve as 
\begin{equation}\label{eq:veloc}
\mathscr{\vv V}(\rho,t) = \vec x_t(\rho,t)\quad 
\text{and}\quad
\mathscr{V}(\rho,t) = \mathscr{\vv V}(\rho,t)\cdot\vec\nu(\rho,t)\qquad\forall\rho\in\mathbb{I},
\end{equation}
respectively.

We have the following lemma for the time derivative of the curvature. 
\begin{lemma}\label{lem:dtkappa} The time derivative of the curvature in \eqref{eq:cur} satisfies 
\begin{equation}\label{eq:dtkappa}
\varkappa_t  = \mathscr{V}_{ss} + \mathscr{V}\varkappa^2 + (\mathscr{\vv V}\cdot\vec\tau)\varkappa_s.
\end{equation}
\end{lemma}
\begin{proof}
The desired result was shown in \cite[Lemma 39(i)]{Barrett20} for an evolving $n$-dimensional hypersurface in $\bR^{n+1}$. For the benefit of reader, here we present a direct proof in the case of an evolving planar curve. 
Taking the time derivative of $\varkappa$ directly gives 
\begin{equation}\label{eq:td1}
\varkappa_t = \bigl(\vec x_{ss}\cdot\vec\nu\bigr)_t = (\vec x_{ss})_t\cdot\vec\nu + \vec x_{ss}\cdot\vec\nu_t=(\vec x_{ss})_t\cdot\vec\nu,
\end{equation}
where we used the fact that $\vec x_{ss}\cdot\vec\nu_t = \varkappa\,\vec\nu\cdot\vec\nu_t = \frac12 \varkappa\,(|\vec\nu|^2)_t = 0$. 

We recall the Frenet--Serret formulas
\begin{equation}\label{eq:dstaun}
\vec\tau_s = \varkappa\,\vec\nu,~\qquad \vec\nu_s = -\varkappa\,\vec\tau.
\end{equation}
Moreover, on recalling \eqref{eqn:taun} and \eqref{eq:veloc} we obtain that
\begin{equation} \label{eq:dtrho}
\partial_t\,\partial_s 
= \left(\frac{\partial_\rho}{|\vec x_\rho|}\right)_t = \frac{\partial_\rho\,\partial_t}{|\vec x_\rho|} + \left(\frac{1}{|\vec x_\rho|}\right)_t \partial_\rho
=\partial_s\,\partial_t - (\mathscr{\vv V}_s\cdot\vec\tau)\partial_s.
\end{equation}

Using \eqref{eq:dtrho} in \eqref{eq:td1} for $\vec\tau_s$ and recalling \eqref{eq:dstaun}, we get 
\begin{equation} \label{eq:dtkappa1}
\varkappa_t = (\vec \tau_{s})_t\cdot\vec\nu  = \bigl[(\vec\tau_t)_s -(\mathscr{\vv V}_s\cdot\vec\tau)\,\vec\tau_s\bigr]\cdot\vec\nu =[\vec\tau_t]_s\cdot\vec\nu - (\mathscr{\vv V}_s\cdot\vec\tau)\,\varkappa.
\end{equation}
Again using \eqref{eq:dtrho} for $\vec x$, we can compute 
\begin{align}
(\vec x_s)_t &= (\vec x_t)_s - (\mathscr{\vv V}_s\cdot\vec\tau)\,\vec\tau=\mathscr{\vv V}_s - (\mathscr{\vv V}_s\cdot\vec\tau)\,\vec\tau\nn\\
&=(\mathscr{\vv V}_s\cdot\vec\nu)\,\vec\nu = (\mathscr{\vv V}\cdot\vec\nu)_s\,\vec\nu - (\mathscr{\vv V}\cdot\vec\nu_s)\,\vec\nu = \left(\mathscr{V}_s + (\mathscr{\vv V}\cdot\vec\tau)\,\varkappa\right)\vec\nu.\label{eq:dttau}
\end{align}

Inserting \eqref{eq:dttau} into \eqref{eq:dtkappa1} and using \eqref{eq:dstaun},  we obtain 
\begin{align*}
\varkappa_t & = \left[\left(\mathscr{V}_s + (\mathscr{\vv V}\cdot\vec\tau)\,\varkappa\right)\,\vec\nu\right]_s\cdot\vec\nu - (\mathscr{\vv V}_s\cdot\vec\tau)\,\varkappa\nn\\
&= \mathscr{V}_{ss} + (\mathscr{\vv V}_s\cdot\vec\tau)\,\varkappa + (\mathscr{\vv V}\cdot\vec\tau_s)\,\varkappa + (\mathscr{\vv V}\cdot\vec\tau)\,\varkappa_s - (\mathscr{\vv V}_s\cdot\vec\tau)\,\varkappa\nn\\
& = \mathscr{V}_{ss}  + (\mathscr{\vv V}\cdot\vec\nu)\,\varkappa^2 + (\mathscr{\vv V}\cdot\vec\tau)\varkappa_s,\nn
\end{align*}
which is exactly \eqref{eq:dtkappa}. 
\end{proof}

In this paper, we consider the $L^2$-gradient flow of the energy 
\begin{equation}\label{eq:energy}
E(\Gamma(t)) = \frac{1}{2}\int_{\Gamma(t)}\varkappa^2\rd s + \lambda\,|\Gamma(t)|,\qquad\lambda\geq 0,
\end{equation}
where $|\Gamma(t)|$ represents the length of the curve. Observe that \eqref{eq:energy} combines the bending energy, involving the curvature of $\Gamma(t)$, with a length contribution that penalizes growth if $\lambda>0$. 
We remark that in \eqref{eq:energy}, and throughout this paper where no
confusion can arise, we use a slight abuse of notation in that we identify
$\varkappa \circ \vec x^{-1}$ with $\varkappa$ on $\Gamma(t) = \vec x(\mathbb I,
t)$.
The desired gradient flow of \eqref{eq:energy} is given by 
\begin{equation}\label{eq:willmore}
\mathscr{V} = -\varkappa_{ss}-\frac{1}{2}\varkappa^3 + \lambda\varkappa,
\end{equation}
see e.g.\ \cite{BGN07}, and obeys the energy dissipation law
\begin{equation}\label{eq:strongenergylaw}
\ddt\, E(\Gamma(t)) =\int_{\Gamma(t)}(\varkappa_{ss}+\frac{1}{2}\varkappa^3-\lambda\varkappa)\,\mathscr{V}\rd s =  -\int_{\Gt}\mathscr{V}^2\rd s\leq 0.
\end{equation}

On recalling Lemma~\ref{lem:dtkappa}, we notice that any parameterization that satisfies \eqref{eq:willmore} also solves the following system
\begin{subequations}\label{eqn:reWillmore}
\begin{align}\label{eq:reWillmore1}
\mathscr{V} &= -\varkappa_{ss} -\frac{1}{2}\varkappa^3 + \lambda\kappa,\\
\label{eq:reWillmore2}
\varkappa_t &= \mathscr{V}_{ss} + \mathscr{V}\,\varkappa^2 + (\vec x_t\cdot\vec\tau)\,\varkappa_s,\\
\label{eq:reWillmore3}
\vec x_t\cdot\vec\nu&=\mathscr{V},\\
\kappa\,\vec\nu &= \vec x_{ss},
\label{eq:reWillmore4}
\end{align}
\end{subequations}
where we introduce an additional curvature variable in \eqref{eq:reWillmore4}. The two curvatures are the same on the continuous level, but will have different approximations on the discrete level.  In order for \eqref{eqn:reWillmore} to be complete, we also need to include the initial parameterization $\vec x(\cdot, 0)$, and the initial curvature $\varkappa(\cdot,0)$, defined via \eqref{eq:cur}.

We recall that the original BGN method from \cite{BGN07} is based on the weak
formulation
\[
\vec x_t\cdot\vec\nu= -\varkappa_{ss} -\frac{1}{2}\varkappa^3 + \lambda\varkappa, \quad
\varkappa\,\vec\nu = \vec x_{ss}.
\]
So compared to this system of two PDEs involving the curvature
$\varkappa$ and the parameterization $\vec x$, the formulation
\eqref{eqn:reWillmore} represents a system of four PDEs with the
additional variables $\mathscr{V}$ and $\kappa$. But crucially, both
formulations do not prescribe a tangential velocity, and both feature the side
constraint \eqref{eq:bgn} that leads to an equidistribution property under
discretization.

\subsection{The weak formulation}\label{sec:weakform}
We multiply \eqref{eq:reWillmore2} with a test function $\chi\in H^1(\mathbb{I})$, integrate over $\Gamma(t)$ and perform integration by parts to get
\begin{align}
&\int_\Gt\varkappa_t\,\chi\rd s =\int_\Gt(\mathscr{V}_{ss} + \mathscr{V}\,\varkappa^2+(\vec x_t\cdot\vec\tau)\,\varkappa_s)\,\chi\rd s \nn \\
&=-\int_\Gt\mathscr{V}_s~\chi_s\rd s + \int_\Gt\mathscr{V}\,\varkappa^2~\chi\rd s+ \int_\Gt(\vec x_t\cdot\vec\tau)\,\varkappa_s~\chi\rd s.\label{eq:imw1}
\end{align}
We can reformulate the last term in \eqref{eq:imw1} as
\begin{align}
&\int_\Gt(\vec x_t\cdot\vec\tau)\,\varkappa_s\,\chi\rd s = \frac{1}{2}\int_\Gt(\vec x_t\cdot\vec\tau)\,[\varkappa_s~\chi - \varkappa ~\chi_s + (\varkappa\,\chi)_s]\rd s \nn\\
& =\frac{1}{2}\int_\Gt(\vec x_t\cdot\vec\tau)\,(\varkappa_s~\chi - \varkappa ~\chi_s)\rd s-\frac{1}{2}\int_\Gt[(\vec x_t)_s\cdot\vec\tau + \vec x_t\cdot\vec\tau_s]\,\varkappa\,\chi\rd s\nn\\
&=\frac{1}{2}\int_\Gt(\vec x_t\cdot\vec\tau)\,(\varkappa_s~\chi - \varkappa ~\chi_s)\rd s-\frac{1}{2}\int_\Gt([\vec x_t]_s\cdot\vec\tau)\,\varkappa\,\chi\rd s \nn \\
&\qquad\qquad - \frac{1}{2}\int_\Gt(\vec x_t\cdot\vec\nu)\,\varkappa^2\,\chi\rd s,
\label{eq:imw2}
\end{align}
where in the second equality we used integration by parts and where the last equality results from the fact that $\vec\tau_s = \varkappa\,\vec\nu$. Combining \eqref{eq:imw1} and \eqref{eq:imw2} leads to 
\begin{align}
&\int_\Gt\varkappa_t\,\chi\rd s = -\int_\Gt\mathscr{V}_s\,\chi_s\rd s + \frac{1}{2}\int_\Gt\mathscr{V}\,\varkappa^2\,\chi\rd s\nn\\ &\qquad -\frac{1}{2}\int_\Gt([\vec x_t]_s\cdot\vec\tau)\,\varkappa\,\chi\rd s+  \frac{1}{2}\int_\Gt(\vec x_t\cdot\vec\tau)\,(\varkappa_s~\chi - \varkappa ~\chi_s)\rd s.\label{eq:imw3}
\end{align}
We note that the last term in \eqref{eq:imw3} is antisymmetric in terms of $\varkappa$ and $\chi$, which will allow us to establish a stability estimate on the discrete level. A similar reformulation for the convective term in a model for two-phase Navier--Stokes flow was utilized in \cite{BGN15stable} to obtain an unconditionally stable unfitted finite element approximation. Later this approach was
adapted by the present authors to a situation with a moving fitted bulk mesh 
in \cite{GNZ24ale}. It turns out that the same ideas will be relevant for the
stable discretization of \eqref{eq:imw3}.

We introduce the $L^2$-inner product over $\mathbb{I}$ as $\bigl(\cdot,\cdot\bigr)$ and then propose the weak formulation for the system \eqref{eqn:reWillmore} as follows.  Given the initial parameterization $\vec x(\cdot,0)$ of
$\Gamma(0) = \vec x(\mathbb I,0)$
with curvature $\varkappa(\cdot,0)$, for each $t\in(0,T]$, we find $\mathscr{V}(t)\in H^1(\mathbb{I})$, $\varkappa(t)\in H^1(\mathbb{I})$, $\vec x(t)\in [H^1(\mathbb{I})]^2$ and $\kappa(t)\in H^1(\mathbb{I})$ such that
\begin{subequations}\label{eqn:weak}
\begin{align}\label{eq:weak1}
&\bigl(\mathscr{V},~\varphi\,|\vec x_\rho|\bigr) = \bigl(\varkappa_\rho,~\varphi_\rho\,|\vec x_\rho|^{-1}\bigr) - \frac{1}{2}\bigl(\varkappa^3,~\varphi\,|\vec x_\rho|\bigr) \\
&\hspace{4cm} + \lambda\bigl(\kappa,~\varphi\,|\vec x_\rho|\bigr)\qquad\forall\varphi\in H^1(\mathbb{I}),\nn\\
\label{eq:weak2}
& \bigl(\varkappa_t,~\chi\,|\vec x_\rho|\bigr) + \frac{1}{2}\bigl([\vec x_t]_\rho\cdot\vec\tau,~\varkappa\,\chi\bigr) = -\bigl(\mathscr{V}_\rho,~\chi_\rho\,|\vec x_\rho|^{-1}\bigr) \\ 
&\hspace{2.2cm} + \frac{1}{2}\bigl(\mathscr{V}\,\varkappa^2,~\chi\,|\vec x_\rho|\bigr) + \frac{1}{2}\bigl(\vec x_t\cdot\vec\tau,~(\varkappa_\rho\,\chi-\varkappa\,\chi_\rho)\bigr) \qquad \forall\chi\in H^1(\mathbb{I}),\nn\\
\label{eq:weak3}
&\bigl(\vec\nu\cdot\vec x_t,~\xi\,|\vec x_\rho|\bigr) =\bigl(\mathscr{V},~\xi\,|\vec x_\rho|\bigr) \qquad\forall\xi\in H^1(\mathbb{I}),\\
&\bigl(\kappa\,\vec\nu,~\vec\eta\,|\vec x_\rho|\bigr) = -\bigl(\vec x_\rho,~\vec\eta_\rho\,|\vec x_\rho|^{-1}\bigr)\qquad\forall\vec\eta\in [H^1(\mathbb{I})]^2.\label{eq:weak4}
\end{align}
\end{subequations}
Here we note that \eqref{eq:weak2} is straightforward in view of \eqref{eq:imw3}. While \eqref{eq:weak1}, \eqref{eq:weak3} and \eqref{eq:weak4}  can be obtained by multiplying \eqref{eq:reWillmore1}, \eqref{eq:reWillmore3} and \eqref{eq:reWillmore4} with suitable test functions, and then integrating over $\Gamma(t)$. We observe that \eqref{eq:weak4} under discretization yields the BGN tangential velocity, which leads to an equidistribution property on the semidiscrete level.

\begin{theorem} The weak solution of \eqref{eqn:weak} satisfies 
\begin{equation*} 
\ddt\left[\dfrac{1}{2}\bigl(\varkappa^2,~\,|\vec x_\rho|\bigr)+ \lambda\,|\Gamma(t)|\right] + \bigl(\mathscr{V},~\mathscr{V}\,|\vec x_\rho|\bigr) =0,
\end{equation*}
which implies the energy dissipation law within the weak formulation. 
\end{theorem}
\begin{proof}
We choose $\varphi=\mathscr{V}$ in \eqref{eq:weak1}, $\chi = \varkappa$ in \eqref{eq:weak2}, $\xi = \lambda\,\kappa$ in \eqref{eq:weak3} and $\vec\eta=\lambda\,\vec x_t$ in \eqref{eq:weak4}, and combine these four equations to obtain that
\begin{align*}
0&=\bigl(\mathscr{V},~\mathscr{V}\,|\vec x_\rho|\bigr)+\lambda\bigl(\vec x_\rho,~[\vec x_t]_\rho\,|\vec x_\rho|^{-1}\bigr) + \bigl(\varkappa_t,~\varkappa\,|\vec x_\rho|\bigr) +\frac{1}{2}\bigl([\vec x_t]_\rho\cdot\vec\tau,~\varkappa^2\bigr)\nn\\
&= \bigl(\mathscr{V},~\mathscr{V}\,|\vec x_\rho|\bigr)+ \lambda\ddt|\Gamma(t)|+ \frac{1}{2}\ddt\bigl(\varkappa^2,~|\vec x_\rho|\bigr),\nn
\end{align*}
as claimed. 
\end{proof}

\section{Semidiscrete scheme}
\label{sec:semi}

We consider a uniform partitioning of the reference domain as $\mathbb{I}=\bigcup_{j=1}^{J}\mathbb{I}_j=\bigcup_{j=1}^J[\rho_{j-1},\rho_{j}]$ with $\rho_j= jh,\, h = 1/J$. Note here that $\rho_0$ and $\rho_J$ are identical.  We then introduce the finite element space
\begin{equation*}
V^h:=\left\{\chi\in C^0(\mathbb{I}): \chi\big|_{\mathbb{I}_j}\quad\mbox{is affine}\quad j = 1,\ldots, J \right\}\subset H^1(\mathbb{I}).\nn
\end{equation*}
We also introduce the mass-lumped $L^2$-inner product $(\cdot,\cdot)^h$ as
\begin{equation}
\bigl( u, v \bigr)^h = \frac{h}{2}\sum_{j=1}^J 
\left[(u\cdot v)(\rho_j^-) + (u\cdot v)(\rho_{j-1}^+)\right],\label{eq:masslumped}
\end{equation}
for scalar or vector valued functions $u,v$, which are piecewise continuous 
with possible jumps at the nodes $\{\rho_j\}_{j=1}^J$, and 
$u(\rho_j^\pm)=\underset{\delta\searrow 0}{\lim}\ u(\rho_j\pm\delta)$. 

Let $(\vec X^h(t))_{t\in[0,T]}$ with $\vec X^h(t)\in [V^h]^2$ be an approximation to $(\vec x(t))_{t\in[0,T]}$. We also define $\Gamma^h(t) = \vec X^h(\mathbb{I},t)$ and assume that 
\begin{equation*}
|\vec X_\rho^h|>0\quad\mbox{in}\quad\mathbb{I}\quad \mbox{for all}\quad t\in[0,T].\nn
\end{equation*}
We also set 
\begin{equation*}
\vec\tau^h = \vec X^h_\rho\,|\vec X_\rho^h|^{-1},\qquad\vec\nu^h = -(\vec\tau^h)^\perp. 
\end{equation*}

The semidiscrete finite element approximation of the weak formulation \eqref{eqn:weak} is given as follows. Given the initial parameterization $\vec x^h(0)\in V^h$ with curvature $\varkappa^h(0)\in V^h$, for each $t\in(0,T]$, we find $\mathscr{V}^h(t)\in V^h$, $\varkappa^h(t)\in V^h$, $\vec X^h(t)\in [V^h]^2$ and $\kappa^h(t)\in V^h$ such that
\begin{subequations}\label{eqn:semid}
\begin{align}\label{eq:semid1}
&\bigl(\mathscr{V}^h,~\varphi^h\,|\vec X^h_\rho|\bigr) = \bigl(\varkappa^h_\rho,~\varphi^h_\rho\,|\vec X^h_\rho|^{-1}\bigr) - \frac{1}{2}\bigl((\varkappa^h)^3,~\varphi^h\,|\vec X^h_\rho|\bigr) \\
&\hspace{4cm} + \lambda\bigl(\kappa^h,~\varphi^h\,|\vec X^h_\rho|\bigr)\qquad\forall\varphi^h\in V^h,\nn\\
\label{eq:semid2}
& \bigl(\varkappa^h_t,~\chi^h\,|\vec X^h_\rho|\bigr) + \frac{1}{2}\bigl([\vec X^h_t]_\rho\cdot\vec\tau^h,~\varkappa^h\,\chi^h\bigr) = -\bigl(\mathscr{V}^h_\rho,~\chi^h_\rho\,|\vec X^h_\rho|^{-1}\bigr) \\
&\quad  + \frac{1}{2}\bigl(\mathscr{V}^h\,(\varkappa^h)^2,~\chi^h\,|\vec X^h_\rho|\bigr) + \frac{1}{2}\bigl(\vec X^h_t\cdot\vec\tau^h,~(\varkappa^h_\rho\,\chi^h-\varkappa^h\,\chi^h_\rho)\bigr) \qquad \forall\chi^h\in V^h,\nn\\
\label{eq:semid3}
&\bigl(\vec\nu^h\cdot\vec X^h_t,~\xi^h\,|\vec X^h_\rho|\bigr)^h =\bigl(\mathscr{V}^h,~\xi^h\,|\vec X^h_\rho|\bigr) \qquad\forall\xi^h\in V^h,\\
&\bigl(\kappa^h\,\vec\nu^h,~\vec\eta^h\,|\vec X^h_\rho|\bigr)^h = -\bigl(\vec X^h_\rho,~\vec\eta^h_\rho\,|\vec X^h_\rho|^{-1}\bigr)\qquad\forall\vec\eta^h\in [V^h]^2.\label{eq:semid4}
\end{align}
\end{subequations}
We observe that we employ mass-lumping in \eqref{eq:semid4} to obtain the
desired BGN tangential motion, and then also in the first term in
\eqref{eq:semid3} to allow an energy estimate if $\lambda>0$.
Indeed, we can prove the following stability estimate for solutions of the
semidiscrete scheme.
\begin{theorem} \label{thm:semides} Let $(\mathscr{V}^h(t), \varkappa^h(t), \vec X^h(t), \kappa^h(t))$ be a solution to the semidiscrete scheme \eqref{eqn:semid}. Then it holds that
\begin{equation}\label{eq:semides}
\ddt\Bigl[\frac{1}{2}\Bigl((\varkappa^h)^2,~|\vec X^h_\rho|\Bigr)+ \lambda\,|\Gamma^h(t)|\Bigr]= -\bigl(\mathscr{V}^h,~\mathscr{V}^h\,|\vec X_\rho^h|\bigr)\leq 0.
\end{equation}
\end{theorem}
\begin{proof}
Choosing $\varphi^h = \mathscr{V}^h(t)$ in \eqref{eq:semid1}, $\chi^h = \varkappa^h(t)$ in \eqref{eq:semid2}, $\xi^h = \lambda\,\kappa^h(t)$ in \eqref{eq:semid3} and $\vec\eta^h = \lambda\,\vec X^h_t(t)$ in \eqref{eq:semid4} and combining these four equations lead to 
\begin{align*}
0&=\bigl(\varkappa^h_t,~\varkappa^h\,|\vec X_\rho^h|\bigr) + \frac{1}{2}\bigl([\vec X_t^h]_\rho\cdot\vec\tau^h,~(\varkappa^h)^2\bigr) + \lambda\bigl(\vec X_\rho^h,~[\vec X_t^h]_\rho\,|\vec X_\rho^h|^{-1}\bigr)\nn\\
&\qquad\qquad  +\bigl(\mathscr{V}^h,~\mathscr{V}^h\,|\vec X_\rho^h|\bigr)\nn\\
&= \ddt\Bigl[\frac{1}{2}\Bigl((\varkappa^h)^2,~|\vec X^h_\rho|\Bigr)+ \lambda\,|\Gamma^h(t)|\Bigr] + \bigl(\mathscr{V}^h,~\mathscr{V}^h\,|\vec X_\rho^h|\bigr),
\end{align*}
which gives \eqref{eq:semides}.
\end{proof}

Denote by
\[\vec a^h_{j-\frac{1}{2}} = \vec X^h(\rho_j, t) -\vec X^h(\rho_{j-1}, t),\qquad j = 1,\ldots, J,\]
the segments of the polygonal curve $\Gamma^h(t)$.
Then we have the following theorem which shows that the vertices on $\Gamma^h(t)$ will be always equidistributed for $t > 0$, provided that they are not locally parallel.
\begin{theorem}[equidistribution] \label{thm:equid} Let $(\mathscr{V}^h(t), \varkappa^h(t), \vec X^h(t), \kappa^h(t))$ be a solution to the semidiscrete scheme \eqref{eqn:semid}. For a fixed time $t\in(0,T]$, it holds for $j=1,\ldots, J$ that
\begin{equation}\label{eq:equid}
|\vec a^h_{j-\frac{1}{2}}| = |\vec a^h_{j+\frac{1}{2}}|\quad\mbox{if}\quad \vec a_{j-\frac{1}{2}} \nparallel \vec a_{j+\frac{1}{2}}.
\end{equation}
\end{theorem}
\begin{proof}
The desired equidistribution property in \eqref{eq:equid} results from \eqref{eq:semid4}, and the proof can be found in \cite[Remark 2.4]{BGN07}.
\end{proof}

\begin{remark}\label{rem:asemi} It is worthwhile to consider an alternative semidiscrete scheme as follows. Given the initial parameterization $\vec x^h(0)\in V^h$ with curvature $\varkappa^h(0)\in V^h$, for each $t\in(0,T]$, we find $\mathscr{V}^h(t)\in V^h$, $\varkappa^h(t)\in V^h$, $\vec X^h(t)\in [V^h]^2$ and $\kappa^h(t)\in V^h$ such that
\begin{subequations}\label{eqn:asemid}
\begin{align}\label{eq:asemid1}
&\bigl(\mathscr{V}^h,~\varphi^h\,|\vec X^h_\rho|\bigr) = \bigl(\varkappa^h_\rho,~\varphi^h_\rho\,|\vec X^h_\rho|^{-1}\bigr) - \frac{1}{2}\bigl((\kappa^h)^2\varkappa^h,~\varphi^h\,|\vec X^h_\rho|\bigr) \\
&\hspace{4cm} + \lambda\bigl(\kappa^h,~\varphi^h\,|\vec X^h_\rho|\bigr)\qquad\forall\varphi^h\in V^h,\nn\\
\label{eq:asemid2}
& \bigl(\varkappa^h_t,~\chi^h\,|\vec X^h_\rho|\bigr) + \frac{1}{2}\bigl([\vec X^h_t]_\rho\cdot\vec\tau^h,~\varkappa^h\,\chi^h\bigr) = -\bigl(\mathscr{V}^h_\rho,~\chi^h_\rho\,|\vec X^h_\rho|^{-1}\bigr) \\
&\quad  + \frac{1}{2}\bigl(\mathscr{V}^h\,(\kappa^h)^2,~\chi^h\,|\vec X^h_\rho|\bigr) + \frac{1}{2}\bigl(\vec X^h_t\cdot\vec\tau^h,~(\varkappa^h_\rho\,\chi^h-\varkappa^h\,\chi^h_\rho)\bigr) \qquad \forall\chi^h\in V^h,\nn
\end{align}
\end{subequations}
together with \eqref{eq:semid3} and \eqref{eq:semid4}. Observe that the only difference is that the square terms $(\varkappa^h)^2$ in \eqref{eq:semid1} and \eqref{eq:semid2} have been replaced by the terms $(\kappa^h)^2$. 
It is easy to see that solutions to this new scheme also satisfy the stability estimate from Theorem~\ref{thm:semides} and the equidistribution property from Theorem~\ref{thm:equid}.
\end{remark}

\section{Fully discrete schemes}\label{sec:fem}

We further divide the time interval uniformly by $[0,T]=\bigcup_{m=1}^{M}[t_{m-1}, t_{m}]$, where $t_m= m\ttau$ and  $\ttau = T/M$ is the time step size.  Let $\vec X^m\in [V^h]^2$ be an approximation to $\vec x(t_m)$ for $0\leq m\leq M$. We define the polygonal curve $\Gamma^m = \vec X^m(\bar{\mathbb{I}})$. Throughout this section we always assume that 
\[|\vec X_\rho^m|>0\quad\mbox{in}\quad\mathbb{I},\quad  m = 0,1,\ldots, M.\]
We also set $\partial_s = |\vec X_\rho^m|^{-1} \partial_\rho$ and
\begin{equation}\label{eq:dtauv}
\vec\tau^m =\vec X_\rho^m\,|\vec X_\rho^m|^{-1},\qquad \vec\nu^m = -(\vec\tau^m)^\perp,
\end{equation}
as the discrete unit tangent and normal vectors of $\Gamma^m$.

For $m=1,\ldots, M$, we introduce the vertex velocity $\mathscr{\vv V}^m\in [V^h]^2$ with 
\begin{equation}\label{eq:dvelocity}
\mathscr{\vv V}^m(\rho_j) = \frac{\vec X^m(\rho_j)-\vec X^{m-1}(\rho_j)}{\ttau},\,\qquad j=1,\ldots, J;
\end{equation}
and define
\begin{equation}
\mathcal{J}^m  = \frac{|\vec X_\rho^{m-1}|}{|\vec X_\rho^m|}.\label{eq:jacob}
\end{equation}

We observe that \eqref{eq:dvelocity} is a consistent temporal approximation
of $\vec X^h_t(\rho_j, t_m)$. Hence it is natural to assume that
$\mathscr{\vv V}^m_s$ remains bounded throughout the evolution. This yields the
following lemma.
\begin{lemma} \label{lem:Jm} 
If $\max_{1 \leq m \leq M} |\mathscr{\vv V}^m_s| \leq C$ for some constant $C$
independent of $\Delta t$, then it holds that
\begin{equation}
\sqrt{\mathcal{J}^m} = 1 - \frac{1}{2}\vec\tau^m\cdot\mathscr{\vv V}^m_\rho\,|\vec X^m_\rho|^{-1}\,\ttau+ O(\ttau^2),\quad m = 1,\ldots, M.\label{eq:Jm}
\end{equation}
\end{lemma}
\begin{proof}
It follows from \eqref{eq:dvelocity} and \eqref{eq:dtauv} that
\begin{align*}
(\mathcal{J}^m)^2 &= \frac{\vec X^{m-1}_\rho\cdot\vec X^{m-1}_\rho}{|\vec X^m_\rho|^2} = \frac{(\vec X^m_\rho - \ttau\,\mathscr{\vv V}_\rho^m)\cdot(\vec X^m_\rho -\ttau\,\mathscr{\vv V}_\rho^m)}{|\vec X^m_\rho|^2} \nn\\
&=  1- 2\ttau\,\vec\tau^m\cdot\mathscr{\vv V}^m_\rho \,|\vec X^m_\rho|^{-1} + \ttau^2\,|\mathscr{\vv V}_s^m|^2.
\end{align*}
Employing the Taylor expansion $(1-\delta)^{\frac{1}{4}} = 1-\frac{1}{4}\delta + O(\delta^2)$ and using our assumption gives the desired result \eqref{eq:Jm}.
\end{proof}

\begin{remark} \label{rem:V}
We will use the property \eqref{eq:Jm} to construct a linear fully discrete
approximation of \eqref{eqn:semid} that is unconditionally stable.
In view of Theorem~\ref{thm:equid}, the derivative $\vec X^h_t$ will have a
jump at time $t=0$ if the initial curve $\Gamma^h(0)$ does not satisfy
\eqref{eq:equid}. Hence for very nonuniform initial data the assumption
of Lemma~\ref{lem:Jm} does not hold, and so the linear scheme
in theory may not be that accurate. However, in practice the scheme works
extremely well also for such initial data. Nevertheless, we will also propose a
nonlinear, and unconditionally stable, scheme that does not rely on the result
from Lemma~\ref{lem:Jm}.
\end{remark}

\subsection{A stable linear scheme}
Let $\mathscr{V}^m, \varkappa^{m}$ and $\kappa^m$
be the numerical approximations of $\mathscr{V}(\cdot,t)$, $\varkappa(\cdot,t)$ and $\kappa(\cdot, t)$ at $t=t_m$.  We propose the following discretized scheme for \eqref{eqn:semid}: Given the initial data $\vec X^0\in [V^h]^2$ and $\varkappa^0\in V^h$, we set $\vec X^{-1}=\vec X^0$, so that $\mathcal{J}^{0} = 1$, and then for $m=0,\ldots, M-1$, we find $\mathscr{V}^{m+1}\in V^h$, $\varkappa^{m+1}\in V^h$, $\vec X^{m+1}\in [V^h]^2$ and $\kappa^{m+1}\in V^h$ such that
\begin{subequations}\label{eqn:fd}
\begin{align}\label{eq:fd1}
&\bigl(\mathscr{V}^{m+1},~\varphi^h\,|\vec X_\rho^m|\bigr)= \bigl(\varkappa_\rho^{m+1},~\varphi_\rho^h\,|\vec X_\rho^m|^{-1}\bigr)\\
&\hspace{1cm}  - \frac{1}{2}\bigl((\varkappa^m)^2\varkappa^{m+1},~\varphi^h\,|\vec X^m_\rho|\bigr)+ \lambda\bigl(\kappa^{m+1},~\varphi^h\,|\vec X^m_\rho|\bigr)\qquad\forall\varphi^h\in V^h,\nn\\
\label{eq:fd2}
&\bigl(\frac{\varkappa^{m+1}-\varkappa^m\,\sqrt{\mathcal{J}^m}}{\ttau},~\chi^h\,|\vec X^m_\rho|\bigr) = -\bigl(\mathscr{V}^{m+1}_\rho,~\chi^h_\rho\,|\vec X^m_\rho|^{-1}\bigr)\\
&\qquad\qquad  + \frac{1}{2}\bigl((\varkappa^m)^2\,\mathscr{V}^{m+1},~\chi^h\,|\vec X^m_\rho|\bigr)\nn\\
&\qquad\qquad  + \frac{1}{2}\bigl(\vec\tau^m\cdot\frac{\vec X^m - \vec X^{m-1}}{\ttau},~[\varkappa^{m+1}_\rho\,\chi^h-\varkappa^{m+1}\,\chi^h_\rho]\bigr)\qquad\forall\chi^h\in V^h,\nn \\
\label{eq:fd3}
&\bigl(\vec\nu^m\cdot\frac{\vec X^{m+1}-\vec X^m}{\ttau},~\xi^h\,|\vec X^m_\rho|\bigr)^h=\bigl(\mathscr{V}^{m+1},~\xi^h\,|\vec X^m_\rho|\bigr) \qquad\forall\xi^h\in V^h,\\
&\bigl(\kappa^{m+1}\,\vec\nu^m,~\vec\eta^h\,|\vec X^m_\rho|\bigr) ^h= -\bigl(\vec X^{m+1}_\rho,~\vec\eta^h_\rho\,|\vec X^m_\rho|^{-1}\bigr)\qquad\forall\vec\eta^h\in [V^h]^2.\label{eq:fd4}
\end{align}
\end{subequations}
Here on recalling \eqref{eq:Jm}, we can recast the first term in \eqref{eq:fd2}  as 
\begin{align*}
&\bigl(\frac{\varkappa^{m+1}-\varkappa^m\sqrt{\mathcal{J}^m}}{\ttau},~\chi^h\,|\vec X^m_\rho|\bigr)\nn\\
&\quad = \bigl( \frac{\varkappa^{m+1}-\varkappa^m}{\ttau},~\chi^h\,|\vec X_\rho^m|\bigr) + \frac{1}{2}\bigl(\vec\tau^m\cdot\vec{\mathscr{V}}_\rho^m,~\varkappa^m\,\chi^h\bigr) +O(\ttau),\nn
\end{align*}
which is hence a consistent temporal discretization of the first two terms in \eqref{eq:semid2}. A similar technique was recently employed by the authors in
\cite{GNZ24ale} in the context of energy stable fully discrete ALE finite element
approximations for two-phase Navier--Stokes flow.

Clearly \eqref{eqn:fd} leads to a linear system of equations at each time level. Note also that when $\lambda=0$ we can decouple the whole system to first solve \eqref{eq:fd1}, \eqref{eq:fd2} for $\mathscr{V}^{m+1}$ and $\varkappa^{m+1}$, and then solve \eqref{eq:fd3}, \eqref{eq:fd4} for $\vec X^{m+1}$ and $\kappa^{m+1}$.

In the following we aim to show that the scheme \eqref{eqn:fd} admits a unique solution, and that the solution satisfies an unconditional stability estimate. 
To this end, we follow \cite{BGN07} to define the vertex normal $\vec\omega^m\in [V^h]^2$ such that
\[
\bigl(\vec\omega^m,~\vec\xi^h\,|\vec X_\rho^m|\bigr)^h = \bigl(\vec\nu^m,~\vec\xi^h\,|\vec X_\rho^m|\bigr)\qquad\forall\vec\xi^h\in [V^h]^2.
\]
It is not difficult to show that
\begin{equation*}
\vec\omega^m(\vec\rho_j) = \frac{-[\vec a^m_{j+\frac{1}{2}} - \vec a^m_{j - \frac{1}{2}}]^\perp}{|\vec a^m_{j+\frac{1}{2}}| + |\vec a^m_{j-\frac{1}{2}}|},\quad j = 1,\ldots, J,
\end{equation*}
where
\begin{equation}
\label{eq:ajm}
\vec a_{j-\frac{1}{2}}^m = \vec X^m(\rho_j) - \vec X^m(\rho_{j-1}),\qquad j = 1,\ldots, J.
\end{equation}
Moreover, it holds that
\begin{equation}\label{eq:weightnor}
\bigl(\chi^h,~\vec\omega^m\cdot\vec\xi^h\,|\vec X_\rho^m|\bigr)^h = \bigl(\chi^h,~\vec\nu^m\cdot\vec\xi^h\,|\vec X_\rho^m|\bigr)^h\quad\forall\chi^h\in V^h,\quad \vec\xi^h\in [V^h]^2.
\end{equation}

\begin{theorem}[existence and uniqueness]  \label{thm:unquesol} 
Assume that 
\begin{enumerate}[label=$(\mathbf{A \arabic*})$, ref=$\mathbf{A \arabic*}$] 
\item \label{assumpI} $\lambda>0$ or  $\vec\omega^m(\rho_j)\neq\vec 0$ for all $j=1,\ldots, J$; 
\item \label{assupII} ${\rm dim\; span}
\bigl\{\vec\omega^{m}(\rho_j)\bigr\}_{j=1}^J=2$.
\end{enumerate}
Then there exists a unique solution $(\mathscr{V}^{m+1}, \varkappa^{m+1}, \vec X^{m+1}, \kappa^{m+1})\in V^h\times V^h\times [V^h]^2\times V^h$ to the linear system \eqref{eqn:fd}.
\end{theorem}
\begin{proof}
Since \eqref{eqn:fd} is linear and the number of unknown equals the number of equations, it suffices to show that the corresponding homogeneous system has only the zero solution. Thus we consider the following homogeneous system for $(\mathscr{V}, \varkappa, \vec X, \kappa)\in V^h\times V^h\times[V^h]^2\times V^h$ such that
\begin{subequations}\label{eqn:homo}
\begin{align}\label{eq:homo1}
&\bigl(\mathscr{V},~\varphi^h\,|\vec X_\rho^m|\bigr)= \bigl(\varkappa_\rho,~\varphi_\rho^h\,|\vec X_\rho^m|^{-1}\bigr) - \frac{1}{2}\bigl((\varkappa^m)^2\varkappa,~\varphi^h\,|\vec X^m_\rho|\bigr)\nn\\
&\hspace{4cm} + \lambda\bigl(\kappa,~\varphi^h\,|\vec X^m_\rho|\bigr)\qquad\forall\varphi^h\in V^h,\\
\label{eq:homo2}
&\bigl(\frac{\varkappa}{\ttau},~\chi^h\,|\vec X^m_\rho|\bigr) = -\bigl(\mathscr{V}_\rho,~\chi^h_\rho\,|\vec X^m_\rho|^{-1}\bigr)  + \frac{1}{2}\bigl((\varkappa^m)^2\,\mathscr{V},~\chi^h\,|\vec X^m_\rho|\bigr)\nn\\
&\qquad\qquad  + \frac{1}{2}\bigl(\vec\tau^m\cdot\frac{\vec X^m - \vec X^{m-1}}{\ttau},~(\varkappa_\rho\,\chi^h-\varkappa\,\chi^h_\rho)\bigr)\qquad\forall\chi^h\in V^h,\\
\label{eq:homo3}
&\bigl(\vec\nu^m\cdot\frac{\vec X}{\ttau},~\xi^h\,|\vec X^m_\rho|\bigr)^h=\bigl(\mathscr{V},~\xi^h\,|\vec X^m_\rho|\bigr) \qquad\forall\xi^h\in V^h,\\
&\bigl(\kappa\,\vec\nu^m,~\vec\eta^h\,|\vec X^m_\rho|\bigr)^h = -\bigl(\vec X_\rho,~\vec\eta^h_\rho\,|\vec X^m_\rho|^{-1}\bigr)\qquad\forall\vec\eta^h\in [V^h]^2.\label{eq:homo4}
\end{align}
\end{subequations}
We choose $\varphi^h = \ttau\,\mathscr{V}$ in \eqref{eq:homo1}, $\chi^h = \ttau\,\varkappa$ in \eqref{eq:homo2}, $\xi^h = \lambda\,\ttau\,\kappa$ in \eqref{eq:homo3} and $\vec\eta^h = \lambda \vec X$ in \eqref{eq:homo4}, and then combine these four equations to obtain 
\[
\ttau\,\bigl(\mathscr{V},~\mathscr{V}\,|\vec X_\rho^m|\bigr) +\bigl(\varkappa,~\varkappa\,|\vec X_\rho^m|\bigr) + \lambda\bigl(\vec X_\rho,~\vec X_\rho\,|\vec X_\rho^m|^{-1}\bigr) = 0. 
\]
This implies immediately that $\mathscr{V}= 0$ and $\varkappa = 0$.
Substituting the former into \eqref{eq:homo3}, we obtain
\begin{equation}\label{eq:xitozero1}
\bigl(\vec\nu^m\cdot\vec X,~\xi^h\,|\vec X^m_\rho|\bigr)^h = 0\qquad\forall\xi^h\in V^h.
\end{equation}
Choosing $\xi^h = \kappa$ in \eqref{eq:xitozero1} and $\vec\eta^h = \vec X$ in \eqref{eq:homo4} gives rise to 
\[\bigl(\vec X_\rho,~\vec X_\rho\,|\vec X_\rho^m|^{-1}\bigr)=0,\]
which implies that $\vec X=\vec X^c$ is constant. Hence \eqref{eq:xitozero1}, together
with \eqref{eq:masslumped} and \eqref{eq:weightnor}, yields that
\begin{equation*}
0=\bigl(\vec X^c\cdot\vec\nu^m,~\xi^h\,|\vec X_\rho^m|\bigr)^h = \bigl(\vec X^c\cdot\vec\omega^m,~\xi^h\,|\vec X_\rho^m|\bigr)^h\qquad\forall\xi^h\in V^h,\nn
\end{equation*}
and so $\vec X^c =\vec 0$ thanks to assumption \ref{assupII}. 
It remains to show that $\kappa = 0$.
If $\lambda > 0$, then choosing $\varphi^h = \kappa$ in \eqref{eq:homo1} yields
that $\kappa = 0$. Otherwise, we have from \eqref{eq:homo4} with
$\vec\eta^h\in [V^h]^2$ defined via
$\vec\eta^h(\rho_j) = \kappa(\rho_j)\vec\omega^m(\rho_j)$, for $j =1,\ldots,
J$, that
\[
0=\bigl(\kappa\,\vec\nu^m,~\vec\eta^h\,|\vec X^m_\rho|\bigr)^h 
= \frac{h}{2}\sum_{j=1}^J|\kappa(\rho_j)\,\vec\omega(\rho_j)|^2(|\vec a^m_{j-\frac{1}{2}}| + |\vec a^m_{j+\frac{1}{2}}|).
\]
This implies that $\kappa(\rho_j) = 0$ for $j=1,\ldots, J$ on recalling assumption \eqref{assumpI}.
Hence we have shown that $(\mathscr{V}, \varkappa, \vec X, \kappa) = (0,0,\vec
0, 0)$ as required.
\end{proof}

\begin{remark} \label{rem:lambda0}
As mentioned earlier, in the case $\lambda=0$ the system \eqref{eqn:fd} decouples into \eqref{eq:fd1}, \eqref{eq:fd2} for $\mathscr{V}^{m+1}$ and $\varkappa^{m+1}$, and \eqref{eq:fd3}, \eqref{eq:fd4} for $\vec X^{m+1}$ and $\kappa^{m+1}$.
The proof of Theorem~\ref{thm:unquesol} shows that in this case, we obtain the
existence of a unique solution to \eqref{eq:fd1}, \eqref{eq:fd2} even without
the two assumptions \ref{assumpI} and \ref{assupII}.
\end{remark}

\begin{theorem}[unconditional stability]\label{thm:energytab} Let $(\mathscr{V}^{m+1}, \varkappa^{m+1}, \vec X^{m+1}, \kappa^{m+1})\in V^h\times V^h\times [V^h]^2\times V^h$  be a solution to the system \eqref{eqn:fd}. Then it holds that
\begin{equation}\label{eq:des}
\mathcal{E}^{m+1}+\ttau\,\bigl(\mathscr{V}^{m+1},~\mathscr{V}^{m+1}\,|\vec X_\rho^m|\bigr)\leq \mathcal{E}^m,
\end{equation}
where $\mathcal{E}^{m} = \dfrac{1}{2}\bigl((\varkappa^{m})^2,\,|\vec X^{m-1}_\rho|\bigr) +\lambda \bigl( |\vec X_\rho^{m}|,~1\bigr).$
\end{theorem}
\begin{proof}
Setting $\varphi^h = \ttau\,\mathscr{V}^{m+1}$ in \eqref{eq:fd1}, $\chi^h = \ttau\,\varkappa^{m+1}$ in \eqref{eq:fd2}, $\xi^h = \lambda\,\ttau\,\kappa^{m+1}$ in \eqref{eq:fd3} and $\vec\eta^h = \lambda\,(\vec X^{m+1}-\vec X^m)$ in \eqref{eq:fd4} and combining these equations gives rise to 
\begin{align}
&\ttau\,\bigl(\mathscr{V}^{m+1},~\mathscr{V}^{m+1}\,|\vec X_\rho^m|\bigr) + \bigl(\varkappa^{m+1}-\varkappa^m\sqrt{\mathcal{J}^m},~\varkappa^{m+1}\,|\vec X_\rho^m|\bigr)\nn\\
&\qquad\qquad + \lambda\bigl(\vec X_\rho^{m+1},~(\vec X^{m+1}-\vec X^m)_\rho\,|\vec X_\rho^m|^{-1}\bigr) =0.\label{eq:des1}
\end{align}

Using the inequality $a(a-b)\geq \frac{1}{2}(a^2-b^2)$, we have
\begin{align}
\bigl(\varkappa^{m+1}-&\varkappa^m\sqrt{\mathcal{J}^m},~\varkappa^{m+1}\,|\vec X_\rho^m|\bigr) \geq \frac{1}{2}\bigl((\varkappa^{m+1})^2 - (\varkappa^m)^2\mathcal{J}^m,~|\vec X_\rho^m|\bigr)\nn\\
&\qquad = \frac{1}{2}\bigl((\varkappa^{m+1})^2,~|\vec X_\rho^m|\bigr)- \frac{1}{2}\bigl((\varkappa^{m})^2,~|\vec X_\rho^{m-1}|\bigr),\label{eq:des2}
\end{align}
where the last equality follows from \eqref{eq:jacob}. 

We next have
\begin{align}
&\bigl(\vec X_\rho^{m+1}, ~(\vec X^{m+1}-\vec X^m)_\rho\,|\vec X_\rho^m|^{-1}\bigr) \geq \frac{1}{2} \bigl(|\vec X_\rho^{m+1}|^2-|\vec X_\rho^m|^2,~|\vec X_\rho^m|^{-1}\bigr) \nn\\
&\qquad\qquad= \frac{1}{2}\bigl(|\vec X_\rho^{m+1}|^2|\vec X_\rho^m|^{-2}-1,~|\vec X_\rho^{m}|\bigr)\nn\\
&\qquad\qquad\geq \bigl(|\vec X_\rho^{m+1}|-|\vec X_\rho^m|,~1\bigr) = |\Gamma^{m+1}|-|\Gamma^m|,\label{eq:des3}
\end{align}
where we used again the fact $a(a-b)\geq \frac{1}{2}(a^2-b^2)$ for the first inequality as well as $\frac{a^2-1}{2}\geq a-1$ for the second inequality.

Now inserting \eqref{eq:des2} and \eqref{eq:des3} into \eqref{eq:des1} yields the desired stability estimate \eqref{eq:des}.
\end{proof}

Analogously to all other existing fully discrete semi-implicit BGN schemes,
e.g.\ from \cite{BGN07,BGN07variational,pwf,Barrett20}, it is not possible to
prove a direct fully discrete analogue of Theorem~\ref{thm:equid}. But clearly any numerical steady state solution of \eqref{eqn:fd} satisfies the equidistribution property \eqref{eq:equid}. Indeed, in practice an asymptotic equidistribution property for the vertices on $\Gamma^m$ can be observed.

\subsection{A stable nonlinear scheme}
We next consider a nonlinear approximation of \eqref{eqn:semid} which also leads to an unconditional stability estimate. Again with the same discrete initial data as before, for $m=0,\ldots, M-1$, we find $\mathscr{V}^{m+1}\in V^h$, $\varkappa^{m+1}\in V^h$, $\vec X^{m+1}\in [V^h]^2$ and $\kappa^{m+1}\in V^h$ such that
\begin{subequations}\label{eqn:nfd}
\begin{align}\label{eq:nfd1}
&\bigl(\mathscr{V}^{m+1},~\varphi^h\,|\vec X_\rho^m|\bigr)= \bigl(\varkappa_\rho^{m+1},~\varphi_\rho^h\,|\vec X_\rho^m|^{-1}\bigr) - \frac{1}{2}\bigl((\varkappa^m)^2\varkappa^{m+1},~\varphi^h\,|\vec X^m_\rho|\bigr) \nn \\
&\hspace{4cm} + \lambda\bigl(\kappa^{m+1},~\varphi^h\,|\vec X^m_\rho|\bigr)\qquad\forall\varphi^h\in V^h,\\
\label{eq:nfd2}
&\bigl(\frac{\varkappa^{m+1}-\varkappa^m}{\ttau},~\chi^h\,|\vec X^m_\rho|\bigr)+\frac{1}{2}\Bigl(\frac{\vec X^{m+1}_\rho-\vec X^m_\rho}{\ttau}\cdot\vec X^{m+1}_\rho\,|\vec X_\rho^m|^{-1},~\varkappa^{m+1}\,\chi^h\Bigr) \nn \\
&\qquad\qquad = -\bigl(\mathscr{V}^{m+1}_\rho,~\chi^h_\rho\,|\vec X^m_\rho|^{-1}\bigr)  + \frac{1}{2}\bigl((\varkappa^m)^2\,\mathscr{V}^{m+1},~\chi^h\,|\vec X^m_\rho|\bigr)\\
&\qquad\qquad\quad  + \frac{1}{2}\bigl(\vec\tau^m\cdot\frac{\vec X^{m+1} -\vec X^{m}}{\ttau},~[\varkappa^{m+1}_\rho\,\chi^h-\varkappa^{m+1}\,\chi^h_\rho]\bigr)\qquad\forall\chi^h\in V^h,\nn\\
\label{eq:nfd3}
&\bigl(\vec\nu^m\cdot\frac{\vec X^{m+1}-\vec X^m}{\ttau},~\xi^h\,|\vec X^m_\rho|\bigr)^h=\bigl(\mathscr{V}^{m+1},~\xi^h\,|\vec X^m_\rho|\bigr) \qquad\forall\xi^h\in V^h,\\
&\bigl(\kappa^{m+1}\,\vec\nu^m,~\vec\eta^h\,|\vec X^m_\rho|\bigr) ^h= -\bigl(\vec X^{m+1}_\rho,~\vec\eta^h_\rho\,|\vec X^m_\rho|^{-1}\bigr)\qquad\forall\vec\eta^h\in [V^h]^2.\label{eq:nfd4}
\end{align}
\end{subequations}
In view of the second and last terms in \eqref{eq:nfd2}, the introduced method \eqref{eqn:nfd} leads to weakly nonlinear systems of polynomial equations at each time level. Due to the presence of \eqref{eq:nfd4}, we retain the nice property of asymptotic equidistribution for the vertices on $\Gamma^m$.    

We have the following theorem for solutions of \eqref{eqn:nfd}, which mimics the energy dissipation law \eqref{eq:strongenergylaw} on the discrete level. 

\begin{theorem}[unconditional stability]\label{thm:energystabnew}
Let $(\mathscr{V}^{m+1}, \varkappa^{m+1}, \vec X^{m+1}, \kappa^{m+1})\in V^h\times V^h\times [V^h]^2\times V^h$  be a solution to the system \eqref{eqn:nfd}. Then it holds that
\begin{equation}\label{eq:ndes}
\mathcal{\bar{E}}^{m+1}+\ttau\,\bigl(\mathscr{V}^{m+1},~\mathscr{V}^{m+1}\,|\vec X_\rho^m|\bigr)\leq \mathcal{\bar{E}}^m,
\end{equation}
where $\mathcal{\bar{E}}^{m} = \dfrac{1}{2}\bigl((\varkappa^{m})^2,\,|\vec X^{m}_\rho|\bigr) + \lambda \bigl( |\vec X_\rho^{m}|,~1\bigr).$
\end{theorem}
\begin{proof}
We choose $\varphi^h = \ttau\,\mathscr{V}^{m+1}$ in \eqref{eq:nfd1}, $\chi^h = \ttau\,\varkappa^{m+1}$ in \eqref{eq:nfd2}, $\xi^h = \lambda\,\ttau\,\kappa^{m+1}$ in \eqref{eq:nfd3} and $\vec\eta^h = \lambda\,(\vec X^{m+1}-\vec X^m)$ in \eqref{eq:nfd4} and combine these equations to obtain  
\begin{align}\label{eq:ndes1}
&\ttau\,\bigl(\mathscr{V}^{m+1},~\mathscr{V}^{m+1}\,|\vec X_\rho^m|\bigr)+ \bigl(\varkappa^{m+1}-\varkappa^m,~\varkappa^{m+1}\,|\vec X_\rho^m|\bigr) \\&\qquad  + \frac{1}{2}\bigl((\vec X^{m+1}-\vec X^m)_\rho\cdot\vec X^{m+1}_\rho\,|\vec X^{m}_\rho|^{-1},~(\varkappa^{m+1})^2\bigr)\nn \\ &\qquad + \lambda\bigl(\vec X_\rho^{m+1},~(\vec X^{m+1}-\vec X^m)_\rho\,|\vec X_\rho^m|^{-1}\bigr) =0.\nn
\end{align}

It is not difficult to show that
\begin{equation} \label{eq:ndes2}
\bigl(\varkappa^{m+1}-\varkappa^m,~\varkappa^{m+1}\,|\vec X_\rho^m|\bigr)\geq \frac{1}{2}\bigl( (\varkappa^{m+1})^2,~|\vec X_\rho^m|\bigr) - \frac{1}{2}\bigl( (\varkappa^{m})^2,~|\vec X_\rho^m|\bigr).
\end{equation}
Moreover, using a similar technique to \eqref{eq:des3}, we have
\begin{align}
&\bigl((\vec X^{m+1}-\vec X^m)_\rho\cdot\vec X^{m+1}_\rho\,|\vec X^{m}_\rho|^{-1},~(\varkappa^{m+1})^2\bigr)\nn\\
 &\quad\geq \frac{1}{2} \bigl(|\vec X_\rho^{m+1}|^2-|\vec X_\rho^m|^2,~|\vec X_\rho^m|^{-1}\,(\varkappa^{m+1})^2\bigr) \nn\\
&\quad= \frac{1}{2}\Bigl(\frac{|\vec X_\rho^{m+1}|^2}{|\vec X_\rho^m|^{2}}-1,~|\vec X_\rho^{m}|\,(\varkappa^{m+1})^2\Bigr)\geq \bigl(|\vec X_\rho^{m+1}|-|\vec X_\rho^m|,~(\varkappa^{m+1})^2\bigr).
\label{eq:ndes3}
\end{align}
Inserting \eqref{eq:ndes2} and \eqref{eq:ndes3} into \eqref{eq:ndes1} and recalling \eqref{eq:des3},  we obtain \eqref{eq:ndes}. 
\end{proof}

We observe that \eqref{eq:ndes} is different from \eqref{eq:des} in terms of the definition of the discrete energy. However, in practice we observe little difference in the behaviour of the solutions of \eqref{eqn:fd} and \eqref{eqn:nfd}.

To solve the nonlinear systems of equations arising at each time level of  \eqref{eqn:nfd}, we propose a Picard-type iteration as follows. For $m\geq 0$, we set $\vec X^{m+1,0} =\vec X^m$. Then, for $\ell \geq 0$, we find 
$\mathscr{V}^{m+1,\ell+1}\in V^h$, $\varkappa^{m+1,\ell+1}\in V^h$, $\vec X^{m+1,\ell+1}\in [V^h]^2$ and $\kappa^{m+1,\ell+1}\in V^h$ such that
\begin{subequations}\label{eqn:pnfd}
\begin{align}\label{eq:pnfd1}
&\bigl(\mathscr{V}^{m+1,\ell+1},~\varphi^h\,|\vec X_\rho^m|\bigr)- \bigl(\varkappa_\rho^{m+1,\ell+1},~\varphi_\rho^h\,|\vec X_\rho^m|^{-1}\bigr) \\ &\qquad\qquad   + \frac{1}{2}\bigl((\varkappa^m)^2\varkappa^{m+1,\ell+1},~\varphi^h\,|\vec X^m_\rho|\bigr) - \lambda\bigl(\kappa^{m+1,\ell+1},~\varphi^h\,|\vec X^m_\rho|\bigr)=0, \nn \\
\label{eq:pnfd2}
&\bigl(\frac{\varkappa^{m+1,\ell+1}-\varkappa^m}{\ttau},~\chi^h\,|\vec X^m_\rho|\bigr)+ \bigl(\mathscr{V}^{m+1,\ell+1}_\rho,~\chi^h_\rho\,|\vec X^m_\rho|^{-1}\bigr)\\
&\qquad\qquad   - \frac{1}{2}\bigl((\varkappa^m)^2\,\mathscr{V}^{m+1,\ell+1},~\chi^h\,|\vec X^m_\rho|\bigr)\nn\\
&\qquad\qquad +\frac{1}{2\,\ttau}\bigl((\vec X^{m+1,\ell}-\vec X^m)_\rho\cdot\vec X^{m+1,\ell}_\rho\,\,|\vec X_\rho^m|^{-1},~\varkappa^{m+1,\ell+1}\,\chi^h\bigr)\nn\\
&\qquad\qquad  - \frac{1}{2}\bigl(\vec\tau^m\cdot\frac{\vec X^{m+1,\ell} - \vec X^{m}}{\ttau},~(\varkappa^{m+1,\ell+1}_\rho\,\chi^h-\varkappa^{m+1,\ell+1}\,\chi^h_\rho)\bigr)=0,\nn\\
\label{eq:pnfd3}
&\bigl(\vec\nu^m\cdot\frac{\vec X^{m+1,\ell+1}-\vec X^m}{\ttau},~\xi^h\,|\vec X^m_\rho|\bigr)^h-\bigl(\mathscr{V}^{m+1,\ell+1},~\xi^h\,|\vec X^m_\rho|\bigr)=0 ,\\
&\bigl(\kappa^{m+1,\ell+1}\,\vec\nu^m,~\vec\eta^h\,|\vec X^m_\rho|\bigr) ^h+\bigl(\vec X^{m+1,\ell+1}_\rho,~\vec\eta^h_\rho\,|\vec X^m_\rho|^{-1}\bigr)=0,\label{eq:pnfd4}
\end{align}
\end{subequations}
for all $(\varphi^h,\chi^h,\xi^h, \vec\eta^h)\in V^h\times V^h\times V^h\times[V^h]^2$.  We then repeat the above iteration until the following condition holds
\begin{equation*}
\max_{1\leq j\leq J}\left\{|\vec X^{m+1,\ell+1}(\rho_j)-\vec X^{m+1,\ell}(\rho_j)|,~ |\varkappa^{m+1,\ell+1}(\rho_j)-\varkappa^{m+1,\ell}(\rho_j)|\right\}\leq {\rm tol},
\end{equation*}
where ${\rm tol}$ is a chosen tolerance.

\subsection{Alternative stable schemes}

In view of the semidiscrete scheme that was mentioned in Remark~\ref{rem:asemi}, it is also worthwhile to consider alternative fully discrete schemes. For example, a new fully discrete linear method can be stated as follows. With the same discrete initial data, for $m=0,\ldots, M-1$, we find $\mathscr{V}^{m+1}\in V^h$, $\varkappa^{m+1}\in V^h$, $\vec X^{m+1}\in [V^h]^2$ and $\kappa^{m+1}\in V^h$ such that
\begin{subequations}\label{eqn:afd}
\begin{align}\label{eq:afd1}
&\bigl(\mathscr{V}^{m+1},~\varphi^h\,|\vec X_\rho^m|\bigr)= \bigl(\varkappa_\rho^{m+1},~\varphi_\rho^h\,|\vec X_\rho^m|^{-1}\bigr) - \frac{1}{2}\bigl((\kappa^m)^2\varkappa^{m+1},~\varphi^h\,|\vec X^m_\rho|\bigr) \nn \\
&\hspace{4cm} + \lambda\bigl(\kappa^{m+1},~\varphi^h\,|\vec X^m_\rho|\bigr)\qquad\forall\varphi^h\in V^h,\\
\label{eq:afd2}
&\bigl(\frac{\varkappa^{m+1}-\varkappa^m\,\sqrt{\mathcal{J}^m}}{\ttau},~\chi^h\,|\vec X^m_\rho|\bigr) = -\bigl(\mathscr{V}^{m+1}_\rho,~\chi^h_\rho\,|\vec X^m_\rho|^{-1}\bigr)\\
&\qquad\qquad  + \frac{1}{2}\bigl((\kappa^m)^2\,\mathscr{V}^{m+1},~\chi^h\,|\vec X^m_\rho|\bigr)\nn\\
&\qquad\qquad  + \frac{1}{2}\bigl(\vec\tau^m\cdot\frac{\vec X^m - \vec X^{m-1}}{\ttau},~[\varkappa^{m+1}_\rho\,\chi^h-\varkappa^{m+1}\,\chi^h_\rho]\bigr)\qquad\forall\chi^h\in V^h,\nn 
\end{align}
\end{subequations}
together with \eqref{eq:fd3} and \eqref{eq:fd4}. It is straightforward to prove existence, uniqueness and unconditional stability for solutions of the new scheme \eqref{eqn:afd} in a similar manner to Theorems~\ref{thm:unquesol} and \ref{thm:energytab}. Moreover, the nonlinear scheme \eqref{eqn:nfd} with $(\varkappa^m)^2$ replaced by $(\kappa^m)^2$ can also be considered, which will once again be
unconditionally stable.

To end this section, we next generalize our approximation to the length-preserving Willmore flow of planar curves. The new flow is  given by 
\begin{equation}
\mathscr{V} = -\varkappa_{ss} -\frac{1}{2}\varkappa^3 + \lambda(t)\varkappa,\label{eq:lengthWillmore} 
\end{equation}
where $\lambda(t)$ is a Lagrange multiplier for the length-preserving constraint 
\[\ddt L(t) = \int_{\Gamma(t)}\mathscr{V}\varkappa\,{\rm d}s = 0.\]
The considered flow can be shown to satisfy 
\begin{equation}    
\tfrac{1}{2}\ddt \int_{\Gamma(t)}\varkappa^2\,{\rm d}s =-\int_{\Gamma(t)}\mathscr{V}^2\,{\rm d}s\leq 0. 
\end{equation}

Based on the method \eqref{eqn:fd}, now the new stable approximation for the geometric flow \eqref{eq:lengthWillmore} can be given as follows. With the same initial discrete data, and then for $m=0,\ldots, M-1$, we find $\mathscr{V}^{m+1}\in V^h$, $\varkappa^{m+1}\in V^h$, $\lambda^{m+1}\in\bR$, $\vec X^{m+1}\in [V^h]^2$ and $\kappa^{m+1}\in V^h$ such that
\begin{align}\label{eqn:lenfd}
&\bigl(\mathscr{V}^{m+1},~\varphi^h\,|\vec X_\rho^m|\bigr)= \bigl(\varkappa_\rho^{m+1},~\varphi_\rho^h\,|\vec X_\rho^m|^{-1}\bigr)\\
&\hspace{1cm}  - \frac{1}{2}\bigl((\varkappa^m)^2\varkappa^{m+1},~\varphi^h\,|\vec X^m_\rho|\bigr)+ \lambda^{m+1}\bigl(\kappa^{m},~\varphi^h\,|\vec X^m_\rho|\bigr)^h\qquad\forall\varphi^h\in V^h,\nn
\end{align}
together with \eqref{eq:fd2}, \eqref{eq:fd3},  \eqref{eq:fd4} and 
\begin{equation}
\bigl(\mathscr{V}^{m+1},~\kappa^m\,|\vec X_\rho^m|\bigr)^h=0.\label{eq:lag}
\end{equation}
It is not difficult to prove existence and uniqueness of a solution for the new scheme \eqref{eqn:lenfd}, \eqref{eq:fd2}, \eqref{eq:fd3}, \eqref{eq:fd4} and \eqref{eq:lag} in a similar manner to Theorem \ref{thm:unquesol}. Moreover, the scheme satisfies an unconditional stability, which is given by 
\[\frac{1}{2}\left((\varkappa^{m+1})^2, |\vec X_\rho^{m}|\right) + \ttau\bigl(\mathscr{V}^{m+1}, \mathscr{V}^{m+1}|\vec X_\rho^m|\bigr)\leq \frac{1}{2}\left((\varkappa^m)^2, |\vec X_\rho^{m-1}|\right).\]

\section{Numerical results}\label{sec:num}

In this section, we present numerical results for  our introduced methods \eqref{eqn:fd} and \eqref{eqn:nfd}. The linear systems arising in \eqref{eqn:fd} and \eqref{eqn:pnfd} are solved with the help of the SparseLU factorizations from the Eigen package \cite{eigenweb}. Throughout the experiments, we start with a nodal interpolation of the initial parameterization $\vec x(\cdot,0)$ to obtain $\vec Y^{0}\in V^h$. We then employ the BGN method \cite{BGN07} with zero normal velocity to compute $\delta\vec Y^0\in [V^h]^2$ and $\kappa^0\in V^h$ such that
\begin{align*}
&\bigl(\vec\nu^{0}_Y\cdot\delta\vec Y^0,~\xi^h\,|\vec Y^0_\rho|\bigr)^h=0 \qquad\forall\xi^h\in V^h,\\
&\bigl(\kappa^{0}\,\vec\nu^0_Y,~\vec\eta^h\,|\vec Y^0_\rho|\bigr) ^h= -\bigl([\vec Y^{0}+\delta\vec Y^0]_\rho,~\vec\eta^h_\rho\,|\vec Y^0_\rho|^{-1}\bigr)\qquad\forall\vec\eta^h\in [V^h]^2, 
\end{align*}
where $\vec\nu^0_Y = -(\vec Y^m_\rho)^\perp\,|\vec Y_\rho^m|^{-1}$. This provides the required initial data $\vec X^0 = \vec Y^0 + \delta\vec Y^0$ and $\varkappa^0 = \kappa^0$. Unless otherwise stated, for the nonlinear method we always choose ${\rm tol} = 10^{-10}$ in our iteration algorithm \eqref{eqn:pnfd}.

\vspace{0.5em}
\noindent
{\bf Example 1}: We begin with a convergence experiment for an expanding circle that was considered in \cite[page 460]{BGN07}. An exact solution to \eqref{eqn:reWillmore} with $\lambda =0 $ is given as
\begin{equation*}
\vec x(\rho,t)= (1+2t)^{\frac{1}{4}}\left(\begin{matrix}
\cos g(\rho)\\
\sin g(\rho)
\end{matrix}\right),\qquad\varkappa(\rho,t) = -(1+2t)^{-\frac{1}{4}}, \quad \rho\in\mathbb{I},\nn
\end{equation*}
for $t\in[0,1]$, where we introduce $g(\rho) = 2\pi\rho + 0.1\sin(2\pi\,\rho)$  to make the initial distribution of nodes nonuniform. We introduce the errors 
\[\norm{\vec X^h - \vec x}_{\infty} = \max_{0\leq m\leq M}\norm{\vec X^m - \vec x(t_m)}_{\infty},\quad \norm{\varkappa^h - \varkappa}_{\infty} = \max_{0\leq m\leq M}\norm{\varkappa^m - \varkappa(t_m)}_{\infty}, \] 
where $\norm{\vec X^m - \vec x(t_m)}_{\infty} = \max_{1 \leq j\leq J}\min_{q\in\mathbb{I}}|\vec X^m(\rho_j) - \vec x(q,t_m)|$ and similarly for $\norm{\varkappa^m - \varkappa(t_m)}_{\infty}$.  We also introduce the error
\[\norm{\kappa^h - \varkappa}_{\infty} = \max_{0\leq m\leq M}\norm{\kappa^m - \varkappa(t_m)}_{\infty}.\]
 The numerical errors and orders of convergence are reported in Table \ref{tb:order}, where we observe a second-order convergence for the two fully discrete schemes in the case of $\ttau = O(h^2)$.  This numerically demonstrates that our introduced schemes are first order in time and second-order in space.

\begin{table}[!htp]
\centering
\def\temptablewidth{0.95\textwidth}
\vspace{-10pt}
\caption{Errors and experimental orders of convergence (EOC) for an expanding circle under Willmore flow with $\lambda =0$ by using the linear scheme \eqref{eqn:fd}(upper panel) and the nonlinear scheme \eqref{eqn:nfd} (lower panel),  where $h_0 = 1/32, \ttau_0=0.04$ and $T=1$.}
{\rule{\temptablewidth}{1pt}}
\begin{tabular*}{\temptablewidth}{@{\extracolsep{\fill}}ccccccc}
$(h,\ \ttau)$  
&$\norm{\vec X^h - \vec x}_{\infty} $ & EOC &$\norm{\varkappa^h - \varkappa}_{\infty}$ & EOC  & $\norm{\kappa^h - \varkappa}_{\infty}$ & EOC  \\[0.2em] \hline
$(h_0, \ttau_0)$
&8.30E-3  &--  &1.38E-2  &--  &4.42E-2 &-- \\[0.2em] \hline 
$(\frac{h_0}{2}, \frac{\ttau_0}{2^2})$
&2.04E-3  &2.02  &3.66E-3 &1.92  &1.11E-2 &1.99 \\[0.2em] \hline 
$(\frac{h_0}{2^2}, \frac{\ttau_0}{2^4})$ 
&5.10E-4  &2.00 &9.27E-4 &1.98 &2.80E-3 &1.99 \\[0.2em]\hline 
$(\frac{h_0}{2^3}, \frac{\ttau_0}{2^6})$ 
&1.27E-4 &2.00 &2.33E-4 &1.99 &7.01E-4 &2.00\\[0.2em]\hline 
$(\frac{h_0}{2^4}, \frac{\ttau_0}{2^8})$ 
&3.18E-5 &2.00 &5.82E-5 &2.00 & 1.75E-4 &2.00 
\end{tabular*}
 {\rule{\temptablewidth}{1pt}}
 {\rule{\temptablewidth}{1pt}}
\begin{tabular*}{\temptablewidth}{@{\extracolsep{\fill}}ccccccc}
$(h_0, \ttau_0)$
&4.38E-3  &--  &4.92E-3  &--  &4.41E-2 &-- \\[0.2em] \hline 
$(\frac{h_0}{2}, \frac{\ttau_0}{2^2})$
&1.09E-3  &2.01  &1.22E-3 &2.01  &1.10E-2 &1.99 \\[0.2em] \hline 
$(\frac{h_0}{2^2}, \frac{\ttau_0}{2^4})$ 
&2.73E-4  &2.00 &3.07E-4 &1.99 &2.80E-3 &1.97 \\[0.2em]\hline 
$(\frac{h_0}{2^3}, \frac{\ttau_0}{2^6})$ 
&6.82E-5 &2.00 &7.67E-5 &1.99 &7.01E-4 &2.00\\[0.2em]\hline 
$(\frac{h_0}{2^4}, \frac{\ttau_0}{2^8})$ 
&1.71E-5 &2.00 &1.92E-5 &2.00 & 1.75E-4 &2.00 
\end{tabular*}
 {\rule{\temptablewidth}{1pt}}
\label{tb:order}
\end{table}

\begin{figure}[!htp]
\centering
\includegraphics[width=0.95\textwidth]{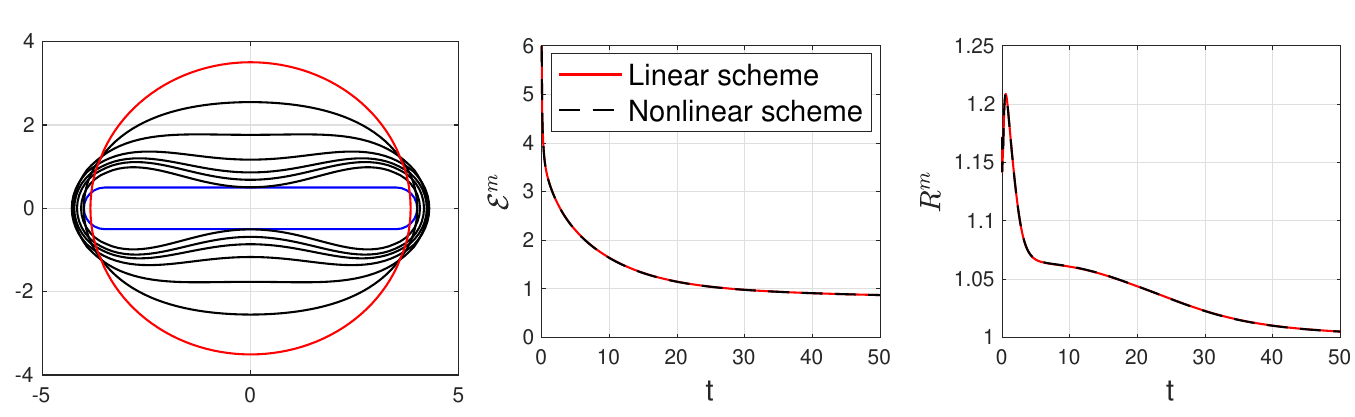}
\caption{[$\lambda =0,~J=128,~\ttau = 10^{-3}$] Evolution for an elongated tube of dimension $8\times 1$,  where we show $\Gamma^m$ at times $t=0,1,2,3,5,10,20,50$ and the time history of the energy $\mathcal{E}^m$ and the mesh ratio $R^m$. }\label{fig:tubelam0}
\end{figure}

\vspace{0.5em}
\noindent
{\bf Example 2}: The next experiment is for an elongated tube of dimension $8\times 1$. We consider the case of $\lambda=0$, and visualize several snapshots of the evolving curve  in Fig.~\ref{fig:tubelam0}. We can see that the curve first become nonconvex and finally evolves towards an expanding circle in order to decrease the energy.  In particular, the discrete energy plots are almost the same for the two schemes and show decreasing trends, which numerically confirms Theorem \ref{thm:energytab} and  Theorem \ref{thm:energystabnew}.

In order to investigate the BGN tangential motion, we further introduce the mesh ratio quantity
\begin{equation*}
{\rm R}^m = \frac{\max_{1\leq j\leq J}|\vec a_{j-\frac{1}{2}}|}{\min_{1\leq j\leq J}|\vec a_{j-\frac{1}{2}}|},\quad m=0,1,\ldots, M, 
\end{equation*}
where $\vec a_{j-\frac{1}{2}}$ is defined as in \eqref{eq:ajm}. The time history of the mesh ratio $R^m$ is plotted in Fig.~\ref{fig:tubelam0} as well for the two schemes. We observe they first increase slightly to 1.2 and then gradually decrease to approximate 1. This implies a uniform distribution of the vertices on the polygonal curve, which substantiates the equidistribution property in Theorem \ref{thm:equid}.
 
\begin{figure}[!htp]
\centering
\includegraphics[width=0.75\textwidth]{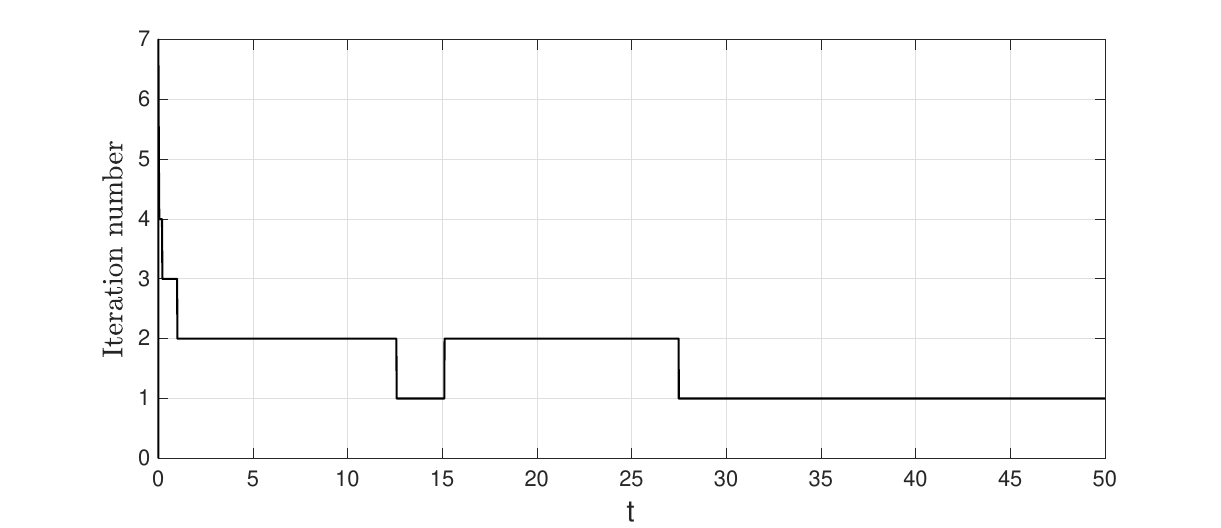}
\caption{The number of iterations for solving the nonlinear scheme \eqref{eqn:nfd} during the experiment in Fig.~\ref{fig:tubelam0}.} \label{fig:picard}
\end{figure}

To test the performance of the Picard iterative algorithm for solving the nonlinear scheme \eqref{eqn:nfd}, we present the number of required iterations in the computation in Fig.~\ref{fig:picard}. Here we observe that apart from the first time steps, the number of iterations is only 2 or 1 due to the weak nonlinearity.

\begin{figure}[!htp]
\centering
\includegraphics[width=0.95\textwidth]{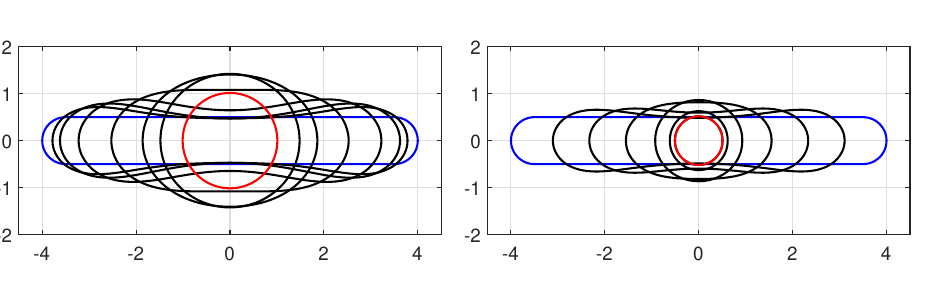}
\includegraphics[width=0.95\textwidth]{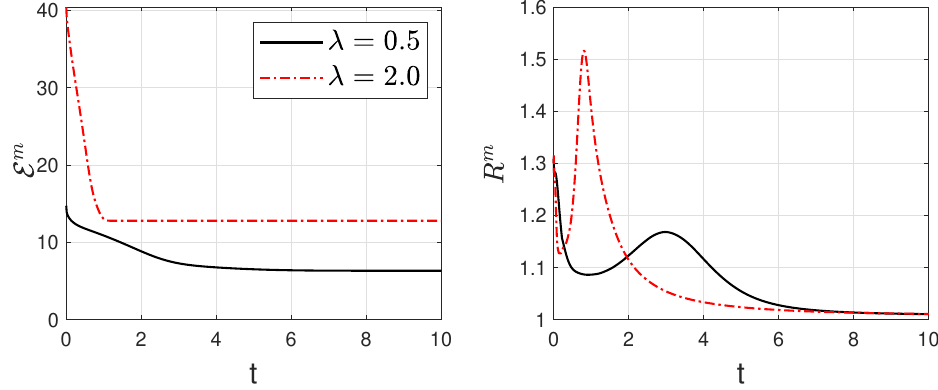}
\caption{[$J=256,~\ttau = 10^{-4}$] Evolution for an elongated tube of dimension $8\times 1$ with $\lambda = 0.5$ (upper left panel) and $\lambda = 2$ (upper right panel), where we visualize $\Gamma^m$ at $t=0,0.2,0.4,1,2,3,4,10$ and $t=0,0.2,0.4,0.6,0.8,1,2,10$, respectively. On the bottom are the plots of the time history of the energy and the mesh ratio $R^m$. }
\label{fig:tubelam}
\end{figure}

As the numerical results from the  linear scheme \eqref{eqn:fd} and the nonlinear scheme \eqref{eqn:nfd} are graphically indistinguishable, in the following we will only report the numerical results from the former.  We next use the same initial setting but consider positive values of $\lambda$ for two new experiments. This length penalization will hinder the curve from growing indefinitely. In Fig.~\ref{fig:tubelam} we visualize $\Gamma^m$ at several times for both $\lambda = 0.5$ and $\lambda=2$. As expected, the curves finally evolve towards steady states which are given by circles of radius 1 and 0.5, respectively.   Moreover, the energy decay and good mesh quality are observed as well in Fig.~\ref{fig:tubelam}.
  
\begin{figure}[!htp]
\centering
\includegraphics[width=0.95\textwidth]{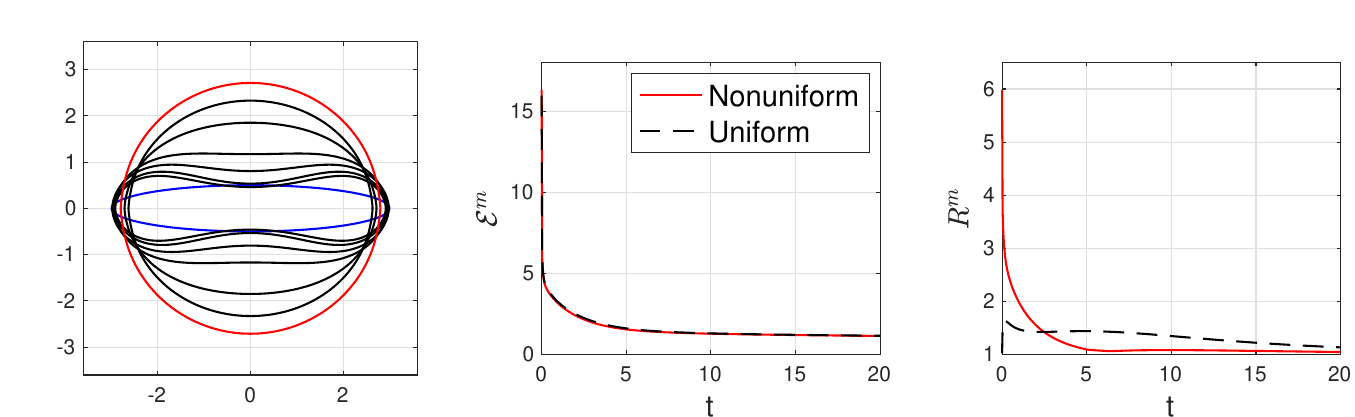}
\caption{[$\lambda = 0,~J = 256,~\ttau = 10^{-3}$] Evolution for an ellipse of dimension $6\times 1$, where we show $\Gamma^m$ at $t=0,0.2,0.4,1,2,5,10,20$, and the plots of the time history of the energy and the mesh ratio $R^m$. For the latter two plots we compare \eqref{eq:ellipse} with a uniform distribution of mesh points on $\Gamma^0$.}
\label{fig:ellipseW}
\end{figure}

\begin{figure}[!htp]
\centering
\includegraphics[width=0.95\textwidth]{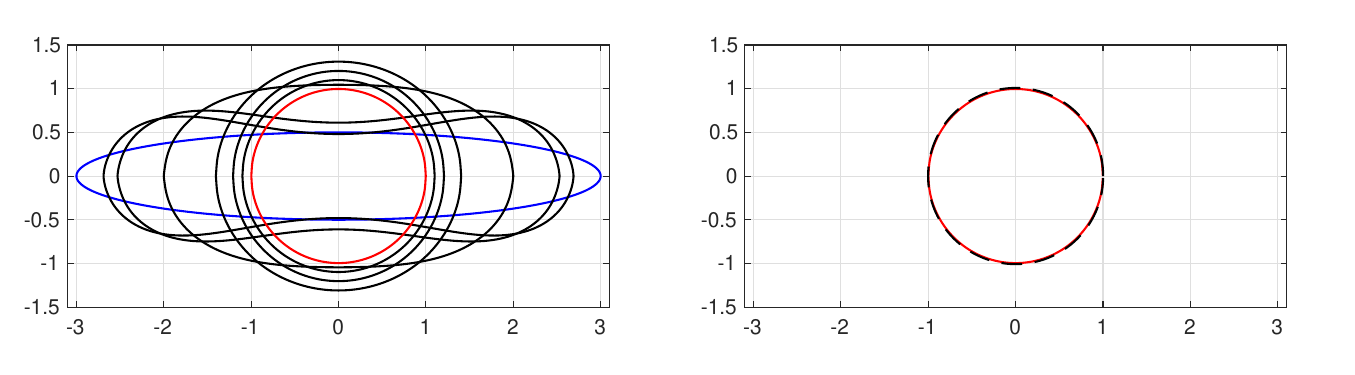}
\includegraphics[width=0.95\textwidth]{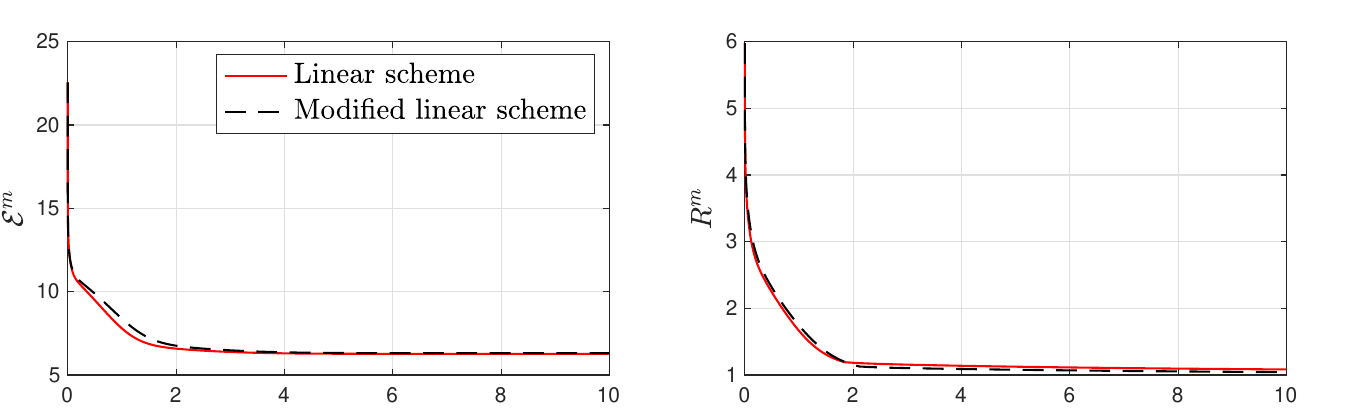}
\caption{[$\lambda = 0.5,~J = 256,~\ttau = 10^{-3}$] Evolution for an ellipse of dimension $6\times 1$ with a nonuniform distribution of vertices on $\Gamma^0$, where we show $\Gamma^m$ at $t=0,0.2,0.4,1,2,3,4,10$. On the upper right are the comparison of $\Gamma^m$ at $t=10$ for the schemes \eqref{eqn:fd} (red solid line) and \eqref{eqn:afd} (black dash line). On the bottom are the plots of the time history of the energy and the mesh ratio $R^m$.}
\label{fig:ellipselam}
\end{figure}

\vspace{0.5em}
\noindent
{\bf Example 3}: In this example, we consider the evolution for an ellipse, which is defined by 
\begin{equation} \label{eq:ellipse}
\vec x(\rho,0) = \left(\begin{array}{ll}
3\cos(2\pi\rho)\\
0.5\sin(2\pi\rho)
\end{array}
\right),\qquad \rho\in\mathbb{I}.
\end{equation}
Notice that a uniform partition of the interval $\mathbb{I}$ then leads to a very nonuniform distribution of vertices on the discrete curve $\Gamma^0$. We conduct an experiment for the case of $\lambda = 0$ and report the numerical results in Fig.~\ref{fig:ellipseW}.  Again we can see that the initially elliptic curve loses its convexity at the beginning, but eventually evolves towards an expanding circle to decrease the energy. Moreover, we observe an asymptotic equidistribution despite the very nonuniform distribution of vertices at the beginning.
  As a comparison, we also perform an experiment where the vertices on the initial $\Gamma^0$ are equally distributed. The resulting curve evolution is indistinguishable from the one shown in Fig.~\ref{fig:ellipseW}, and we observe that the discrete energy plots are almost the same for the two experiments. This demonstrate the robustness of our introduced scheme.

\begin{figure}[!htp]
\includegraphics[width=0.95\textwidth]{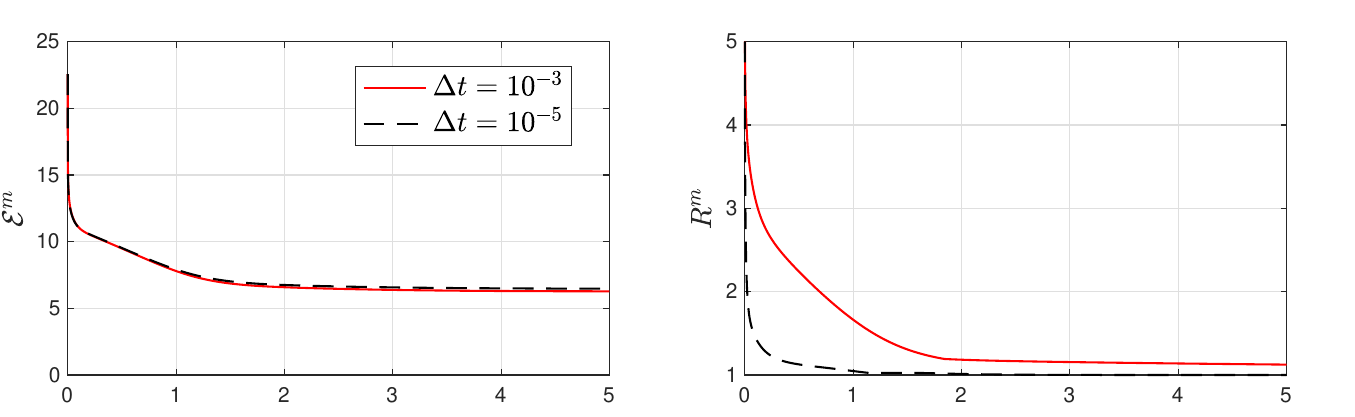}
\caption{[$\lambda = 0.5,~J = 256$] The time history of the energy for the evolution of an ellipse of dimension $6\times 1$ with a nonuniform distribution of vertices on $\Gamma^0$ and different time step sizes.}
\label{fig:smalldt}
\end{figure}

We also conduct an experiment with a very small time step size $\ttau = 10^{-5}$, and compare the results with those obtained from $\ttau = 10^{-3}$. As shown in Fig.~\ref{fig:smalldt}, the discrete energy plots for both experiments are nearly identical, whereas the mesh ratio function converges to 1 significantly faster when using the smaller time step size.

\vspace{0.5em}
\noindent
{\bf Example 4}: In this example, we aim to compare the performance between the linear scheme \eqref{eqn:fd} and its modified variant \eqref{eqn:afd}. We again use the ellipse \eqref{eq:ellipse} for the discrete initial data. We then conduct experiments for the two linear schemes in the case of $\lambda=0.5$. The numerical results are reported in Fig.~\ref{fig:ellipselam}, where we observe that in both cases the ellipse evolves towards a unit circle as the steady state. Moreover, the plots for the monotonically decreasing energy and the mesh ratio are nearly identical. 

\begin{figure}[!htp]
\centering
\includegraphics[width=0.95\textwidth]{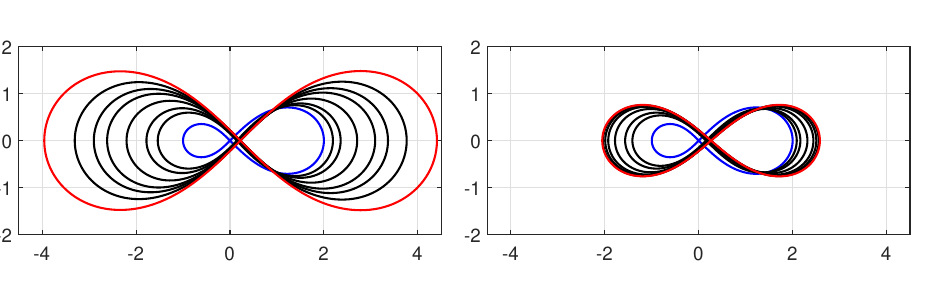}
\includegraphics[width=0.95\textwidth]{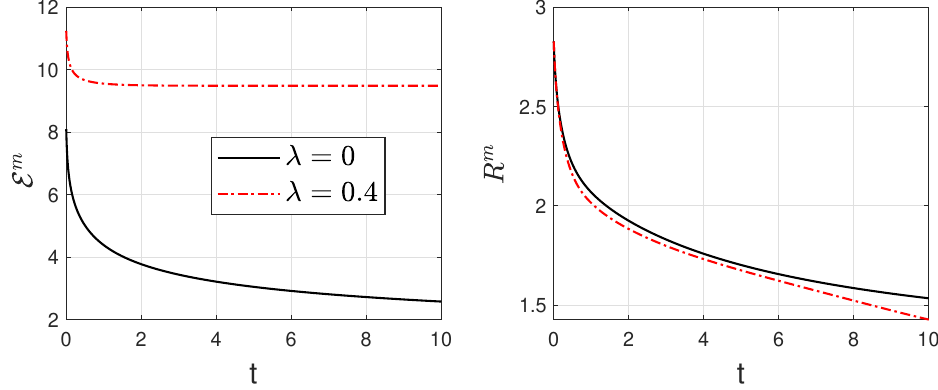}
\includegraphics[width=0.7\textwidth]{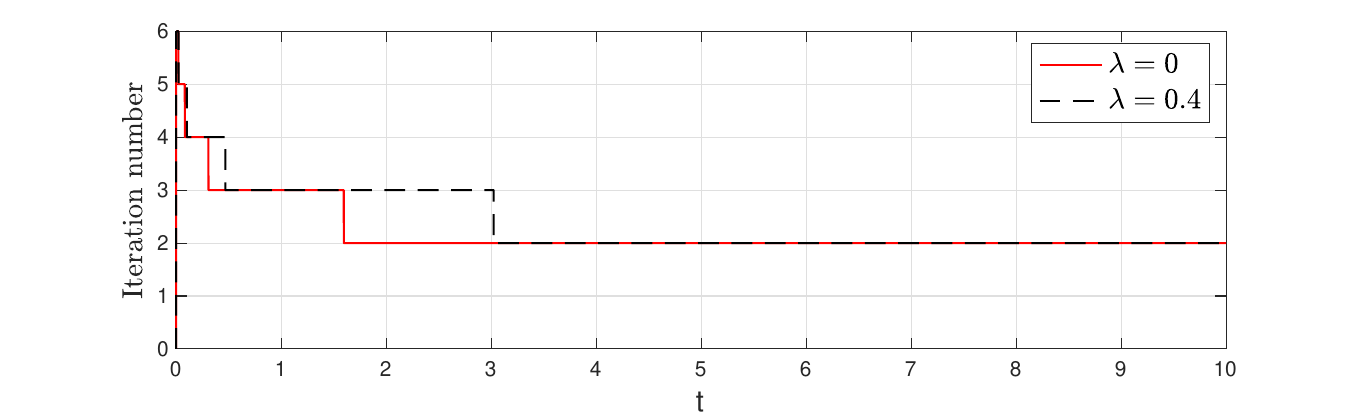}
\caption{[$J=256,~\ttau = 10^{-3}$] Evolution for a 2:1 lemniscate with $\lambda = 0$ (upper left panel) and $\lambda=0.4$ (upper right panel), where we show $\Gamma^m$ at times $t=0,0.2,0.4,1,2,3,5,10$. On the middle are the plots of the time history of the energy and the mesh ratio $R^m$. On the bottom are the iteration numbers when using the nonlinear scheme \eqref{eqn:nfd}.}
\label{fig:lemni}
\end{figure}

\vspace{0.5em}
\noindent
{\bf Example 5}: In our next example, we demonstrate the evolution of a 2:1 lemniscate, which is given by 
\begin{equation*}
\vec x(\rho, 0) = \Bigl(\frac{a(\rho)\,\cos(2\pi\rho)}{1+\sin^2(2\pi\rho)}, ~\frac{a(\rho)\,\cos(2\pi\rho)\sin(2\pi\rho)}{1+\sin^2(2\pi\rho)}\Bigr)^T,\quad\rho\in\mathbb{I},\nn
\end{equation*}
with $a(\rho) = \left\{\begin{array}{ll}
1\quad\mbox{if}\quad \rho\in[\frac{1}{4},\frac{3}{4}],\\
2\quad\mbox{otherwise}.
\end{array}
\right.$ We show the evolution of the lemniscate with $\lambda=0$ and $\lambda = 0.4$ in Fig.~\ref{fig:lemni}. We can see that the nonsymmetric lemniscate evolves to form a symmetric shape. In the case $\lambda=0$ the evolution
approaches a self-similarly expanding symmetric lemniscate, while a steady state can be achieved when a length penalization is applied. Nevertheless, in both cases energy decay and an asymptotic equidistribution property can be observed. This further confirms the excellent performance of our introduced schemes.  

We also examine the performance of the nonlinear scheme \eqref{eqn:nfd} for this  curve evolution. Here we observe for both cases that the required iteration number is about 2 at most times, which is very efficient.

\begin{figure}[!htp]
\centering
\includegraphics[width=0.7\textwidth]{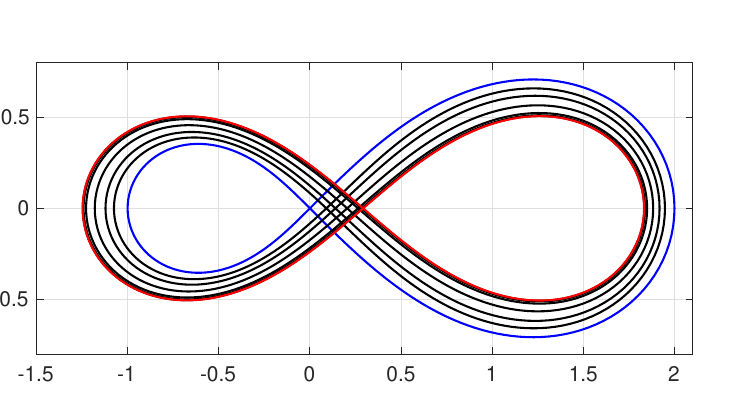}
\includegraphics[width=0.95\textwidth]{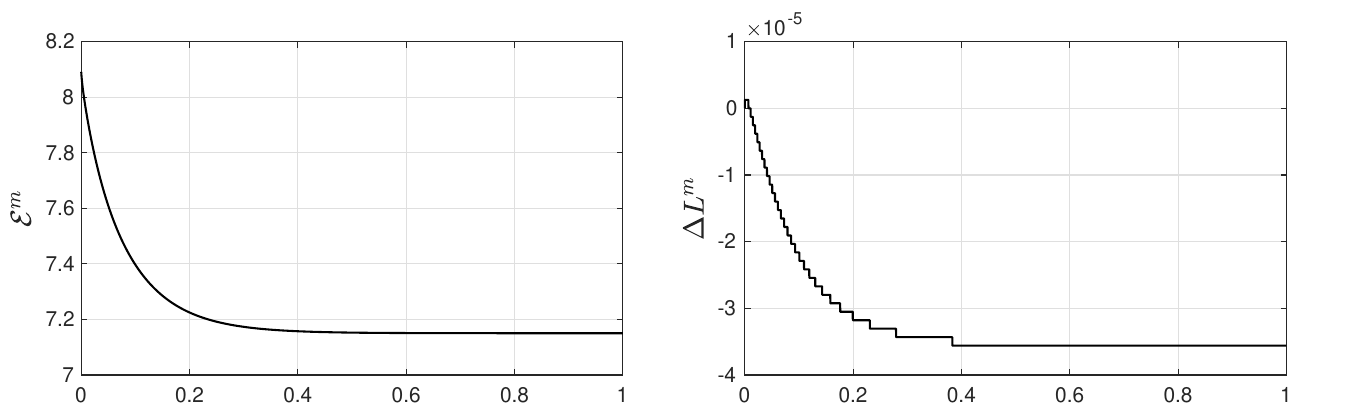}
\caption{[$J=256,~\ttau = 10^{-3}$] Evolution for a 2:1 lemniscate for length-preserving Willmore flow \eqref{eq:lengthWillmore}, where we show $\Gamma^m$ at times $t=0,0.05,0.1,0.2,0.4,0.6,0.8,1$. On the bottom are the plots of the time history of the energy and the relative length change $\Delta L^m$. }
\label{fig:lengthW}
\end{figure}

\vspace{0.5em}
\noindent
{\bf Example 6}: In our final example, we apply the linear scheme \eqref{eqn:lenfd} to the length-preserving Willmore flow of planar curves. The remaining settings are the same as in Example 5. The results are presented in Fig.~\ref{fig:lengthW}, where we observe the expected energy decay of the numerical solution. Moreover, the relative length change, defined by  
\[\Delta L^m = \frac{|\Gamma^m|-|\Gamma^0|}{|\Gamma^0|},\] remains as  small as $10^{-5}$, which indicates that the lengths of the discrete polygonal curves are well preserved.

\section{Conclusions}
\label{sec:con}

We proposed and analyzed two energy-stable fully discrete schemes for Willmore flow of planar curves. The introduced schemes are based on a novel weak formulation that involves an evolution equation for the curvature together with the BGN formulation of the curvature identity. We employed piecewise linear elements in space to obtain a semidiscrete scheme, which satisfies a stability bound and an equidistribution property. Further temporal discretization leads to two distinct fully discrete schemes. The first scheme is linear and unconditionally stable in terms of a discrete energy defined via geometric quantities at two time levels. The second scheme is mildly nonlinear but enjoys an unconditional energy stability where the discrete energy is defined through quantities at a single time level.  We demonstrate a variety of numerical examples to confirm the numerical accuracy and strong performance of the two introduced schemes. 

In the current work, we have only considered the case of planar curves. Extending the presented approach to Willmore flow of surfaces, both in the axisymmetric case and in the general three-dimensional case, are part of our ongoing research.

\section*{Acknowledgement}
This work was partially supported by the National Natural Science Foundation of China No.~12401572 (Q.Z).

\section*{Data availibility}The datasets generated and/or analysed during the current study are available from the corresponding author on reasonable request.
\section*{Declarations}The authors certify that they have no affiliations with or involvement in any organization or entity with any financial interest or non-financial interest in the subject matter or materials discussed in this manuscript.

\bibliographystyle{model1b-num-names}
\bibliography{bib}

\begin{thebibliography}{35}
\expandafter\ifx\csname natexlab\endcsname\relax\def\natexlab#1{#1}\fi
\providecommand{\bibinfo}[2]{#2}
\ifx\xfnm\relax \def\xfnm[#1]{\unskip,\space#1}\fi
\bibitem[{Balzani and Rumpf(2012)}]{BalzaniR12}
\bibinfo{author}{N.~Balzani}, \bibinfo{author}{M.~Rumpf}, \bibinfo{title}{A
  nested variational time discretization for parametric {W}illmore flow},
  \bibinfo{journal}{Interfaces Free Bound.} \bibinfo{volume}{14}
  (\bibinfo{year}{2012}) \bibinfo{pages}{431--454}.
\bibitem[{Bao and Li(2025)}]{BaoL25}
\bibinfo{author}{W.~Bao}, \bibinfo{author}{Y.~Li}, \bibinfo{title}{An
  energy-stable parametric finite element method for the planar {W}illmore
  flow}, \bibinfo{journal}{SIAM J. Numer. Anal.} \bibinfo{volume}{63}
  (\bibinfo{year}{2025}) \bibinfo{pages}{103--121}.
\bibitem[{Barrett et~al.(2007{\natexlab{a}})Barrett, Garcke and
  N\"urnberg}]{BGN07variational}
\bibinfo{author}{J.W. Barrett}, \bibinfo{author}{H.~Garcke},
  \bibinfo{author}{R.~N\"urnberg}, \bibinfo{title}{On the variational
  approximation of combined second and fourth order geometric evolution
  equations}, \bibinfo{journal}{SIAM J. Sci. Comput.} \bibinfo{volume}{29}
  (\bibinfo{year}{2007}{\natexlab{a}}) \bibinfo{pages}{1006--1041}.
\bibitem[{Barrett et~al.(2007{\natexlab{b}})Barrett, Garcke and
  N\"urnberg}]{BGN07}
\bibinfo{author}{J.W. Barrett}, \bibinfo{author}{H.~Garcke},
  \bibinfo{author}{R.~N\"urnberg}, \bibinfo{title}{A parametric finite element
  method for fourth order geometric evolution equations}, \bibinfo{journal}{J.
  Comput. Phys.} \bibinfo{volume}{222} (\bibinfo{year}{2007}{\natexlab{b}})
  \bibinfo{pages}{441--467}.
\bibitem[{Barrett et~al.(2008)Barrett, Garcke and N\"urnberg}]{BGN08willmore}
\bibinfo{author}{J.W. Barrett}, \bibinfo{author}{H.~Garcke},
  \bibinfo{author}{R.~N\"urnberg}, \bibinfo{title}{Parametric approximation of
  {W}illmore flow and related geometric evolution equations},
  \bibinfo{journal}{SIAM J. Sci. Comput.} \bibinfo{volume}{31}
  (\bibinfo{year}{2008}) \bibinfo{pages}{225--253}.
\bibitem[{Barrett et~al.(2010)Barrett, Garcke and N\"urnberg}]{curves3d}
\bibinfo{author}{J.W. Barrett}, \bibinfo{author}{H.~Garcke},
  \bibinfo{author}{R.~N\"urnberg}, \bibinfo{title}{Numerical approximation of
  gradient flows for closed curves in {${\mathbb R}^d$}}, \bibinfo{journal}{IMA
  J. Numer. Anal.} \bibinfo{volume}{30} (\bibinfo{year}{2010})
  \bibinfo{pages}{4--60}.
\bibitem[{Barrett et~al.(2012)Barrett, Garcke and N\"urnberg}]{pwf}
\bibinfo{author}{J.W. Barrett}, \bibinfo{author}{H.~Garcke},
  \bibinfo{author}{R.~N\"urnberg}, \bibinfo{title}{Parametric approximation of
  isotropic and anisotropic elastic flow for closed and open curves},
  \bibinfo{journal}{Numer. Math.} \bibinfo{volume}{120} (\bibinfo{year}{2012})
  \bibinfo{pages}{489--542}.
\bibitem[{Barrett et~al.(2015)Barrett, Garcke and N\"urnberg}]{BGN15stable}
\bibinfo{author}{J.W. Barrett}, \bibinfo{author}{H.~Garcke},
  \bibinfo{author}{R.~N\"urnberg}, \bibinfo{title}{A stable parametric finite
  element discretization of two-phase {Navier--Stokes} flow},
  \bibinfo{journal}{J. Sci. Comput.} \bibinfo{volume}{63}
  (\bibinfo{year}{2015}) \bibinfo{pages}{78--117}.
\bibitem[{Barrett et~al.(2020)Barrett, Garcke and N\"urnberg}]{Barrett20}
\bibinfo{author}{J.W. Barrett}, \bibinfo{author}{H.~Garcke},
  \bibinfo{author}{R.~N\"urnberg}, \bibinfo{title}{Parametric finite element
  approximations of curvature driven interface evolutions},
  \bibinfo{journal}{Handb. Numer. Anal. (Andrea Bonito and Ricardo H. Nochetto,
  eds.)} \bibinfo{volume}{21} (\bibinfo{year}{2020}) \bibinfo{pages}{275--423}.
\bibitem[{Bartels(2013)}]{Bartels13a}
\bibinfo{author}{S.~Bartels}, \bibinfo{title}{A simple scheme for the
  approximation of the elastic flow of inextensible curves},
  \bibinfo{journal}{IMA J. Numer. Anal.} \bibinfo{volume}{33}
  (\bibinfo{year}{2013}) \bibinfo{pages}{1115--1125}.
\bibitem[{Bondarava(2015)}]{Bondarava15}
\bibinfo{author}{A.~Bondarava}, \bibinfo{title}{Stability and error analysis
  for a numerical scheme to approximate elastic flow}, Ph.D. thesis, University
  Magdeburg, \bibinfo{address}{Magdeburg}, \bibinfo{year}{2015}.
\bibitem[{Bretin et~al.(2011)Bretin, Lachaud and Oudet}]{BretinLO11}
\bibinfo{author}{E.~Bretin}, \bibinfo{author}{J.O. Lachaud},
  \bibinfo{author}{E.~Oudet}, \bibinfo{title}{Regularization of discrete
  contour by {W}illmore energy}, \bibinfo{journal}{J. Math. Imaging Vision}
  \bibinfo{volume}{40} (\bibinfo{year}{2011}) \bibinfo{pages}{214--229}.
\bibitem[{Canham(1970)}]{Canham1970minimum}
\bibinfo{author}{P.B. Canham}, \bibinfo{title}{The minimum energy of bending as
  a possible explanation of the biconcave shape of the human red blood cell},
  \bibinfo{journal}{J. Theor. Bio.} \bibinfo{volume}{26} (\bibinfo{year}{1970})
  \bibinfo{pages}{61--81}.
\bibitem[{Chen et~al.(2018)Chen, Lowengrub, Shen, Wang and Wise}]{ChenLSWW18}
\bibinfo{author}{Y.~Chen}, \bibinfo{author}{J.~Lowengrub},
  \bibinfo{author}{J.~Shen}, \bibinfo{author}{C.~Wang},
  \bibinfo{author}{S.~Wise}, \bibinfo{title}{Efficient energy stable schemes
  for isotropic and strongly anisotropic {C}ahn--{H}illiard systems with the
  {W}illmore regularization}, \bibinfo{journal}{J. Comput. Phys.}
  \bibinfo{volume}{365} (\bibinfo{year}{2018}) \bibinfo{pages}{56--73}.
\bibitem[{Clarenz et~al.(2004)Clarenz, Diewald, Dziuk, Rumpf and
  Rusu}]{ClarenzDDRR04}
\bibinfo{author}{U.~Clarenz}, \bibinfo{author}{U.~Diewald},
  \bibinfo{author}{G.~Dziuk}, \bibinfo{author}{M.~Rumpf},
  \bibinfo{author}{R.~Rusu}, \bibinfo{title}{A finite element method for
  surface restoration with smooth boundary conditions},
  \bibinfo{journal}{Comput. Aided Geom. Design} \bibinfo{volume}{21}
  (\bibinfo{year}{2004}) \bibinfo{pages}{427--445}.
\bibitem[{Deckelnick and Dziuk(2009)}]{DeckelnickD09}
\bibinfo{author}{K.~Deckelnick}, \bibinfo{author}{G.~Dziuk},
  \bibinfo{title}{Error analysis for the elastic flow of parametrized curves},
  \bibinfo{journal}{Math. Comp.} \bibinfo{volume}{78} (\bibinfo{year}{2009})
  \bibinfo{pages}{645--671}.
\bibitem[{Deckelnick et~al.(2005)Deckelnick, Dziuk and
  Elliott}]{DeckelnickDE05}
\bibinfo{author}{K.~Deckelnick}, \bibinfo{author}{G.~Dziuk},
  \bibinfo{author}{C.M. Elliott}, \bibinfo{title}{Computation of geometric
  partial differential equations and mean curvature flow},
  \bibinfo{journal}{Acta Numer.} \bibinfo{volume}{14} (\bibinfo{year}{2005})
  \bibinfo{pages}{139--232}.
\bibitem[{Deckelnick and N\"urnberg(2024)}]{cd}
\bibinfo{author}{K.~Deckelnick}, \bibinfo{author}{R.~N\"urnberg},
  \bibinfo{title}{Finite element schemes with tangential motion for fourth
  order geometric curve evolutions in arbitrary codimension},
  \bibinfo{howpublished}{arXiv:2402.16799}, \bibinfo{year}{2024}.
\bibitem[{Desbrun et~al.(1999)Desbrun, Meyer, Schr{\"o}der and
  Barr}]{Desbrun99implicit}
\bibinfo{author}{M.~Desbrun}, \bibinfo{author}{M.~Meyer},
  \bibinfo{author}{P.~Schr{\"o}der}, \bibinfo{author}{A.H. Barr},
  \bibinfo{title}{Implicit fairing of irregular meshes using diffusion and
  curvature flow}, in: \bibinfo{booktitle}{Proceedings of the 26th annual
  conference on Computer graphics and interactive techniques}, pp.
  \bibinfo{pages}{317--324}.
\bibitem[{Duan et~al.(2021)Duan, Li and Zhang}]{Duan2021high}
\bibinfo{author}{B.~Duan}, \bibinfo{author}{B.~Li}, \bibinfo{author}{Z.~Zhang},
  \bibinfo{title}{High-order fully discrete energy diminishing evolving surface
  finite element methods for a class of geometric curvature flows},
  \bibinfo{journal}{Ann. Appl. Math} \bibinfo{volume}{37}
  (\bibinfo{year}{2021}) \bibinfo{pages}{405--436}.
\bibitem[{Dziuk(2008)}]{Dziuk08}
\bibinfo{author}{G.~Dziuk}, \bibinfo{title}{Computational parametric {W}illmore
  flow}, \bibinfo{journal}{Numer. Math.} \bibinfo{volume}{111}
  (\bibinfo{year}{2008}) \bibinfo{pages}{55--80}.
\bibitem[{Dziuk et~al.(2002)Dziuk, Kuwert and Sch{\"a}tzle}]{DziukKS02}
\bibinfo{author}{G.~Dziuk}, \bibinfo{author}{E.~Kuwert},
  \bibinfo{author}{R.~Sch{\"a}tzle}, \bibinfo{title}{Evolution of elastic
  curves in {${\mathbb R}^n$}: {E}xistence and computation},
  \bibinfo{journal}{SIAM J. Math. Anal.} \bibinfo{volume}{33}
  (\bibinfo{year}{2002}) \bibinfo{pages}{1228--1245}.
\bibitem[{Garcke et~al.(2024)Garcke, N\"urnberg and Zhao}]{GNZ24ale}
\bibinfo{author}{H.~Garcke}, \bibinfo{author}{R.~N\"urnberg},
  \bibinfo{author}{Q.~Zhao}, \bibinfo{title}{Arbitrary {L}agrangian--{E}ulerian
  finite element approximations for axisymmetric two-phase flow},
  \bibinfo{journal}{Comput. Math. Appl.} \bibinfo{volume}{155}
  (\bibinfo{year}{2024}) \bibinfo{pages}{209--223}.
\bibitem[{Grinspun(2008)}]{Grinspun08}
\bibinfo{author}{E.~Grinspun}, \bibinfo{title}{A discrete model of thin
  shells}, in: \bibinfo{booktitle}{Discrete differential geometry},
  volume~\bibinfo{volume}{38} of \textit{\bibinfo{series}{Oberwolfach Semin.}},
  \bibinfo{publisher}{Birkh\"auser, Basel}, \bibinfo{year}{2008}, pp.
  \bibinfo{pages}{325--337}.
\bibitem[{Gruber and Aulisa(2020)}]{GruberA20}
\bibinfo{author}{A.~Gruber}, \bibinfo{author}{E.~Aulisa},
  \bibinfo{title}{Computational p-{W}illmore flow with conformal penalty},
  \bibinfo{journal}{ACM Trans. Graph.} \bibinfo{volume}{39}
  (\bibinfo{year}{2020}) \bibinfo{pages}{161}.
\bibitem[{Guennebaud et~al.(2010)Guennebaud, Jacob et~al.}]{eigenweb}
\bibinfo{author}{G.~Guennebaud}, \bibinfo{author}{B.~Jacob}, et~al.,
  \bibinfo{title}{Eigen v3},
  \bibinfo{howpublished}{http://eigen.tuxfamily.org}, \bibinfo{year}{2010}.
\bibitem[{Helfrich(1973)}]{Helfrich73elastic}
\bibinfo{author}{W.~Helfrich}, \bibinfo{title}{Elastic properties of lipid
  bilayers: theory and possible experiments}, \bibinfo{journal}{Z. Naturforsch
  C} \bibinfo{volume}{28} (\bibinfo{year}{1973}) \bibinfo{pages}{693--703}.
\bibitem[{Kov\'{a}cs et~al.(2019)Kov\'{a}cs, Li and Lubich}]{KovacsLL19}
\bibinfo{author}{B.~Kov\'{a}cs}, \bibinfo{author}{B.~Li},
  \bibinfo{author}{C.~Lubich}, \bibinfo{title}{A convergent evolving finite
  element algorithm for mean curvature flow of closed surfaces},
  \bibinfo{journal}{Numer. Math.} \bibinfo{volume}{143} (\bibinfo{year}{2019})
  \bibinfo{pages}{797--853}.
\bibitem[{Kov\'{a}cs et~al.(2021)Kov\'{a}cs, Li and Lubich}]{KovacsLL21}
\bibinfo{author}{B.~Kov\'{a}cs}, \bibinfo{author}{B.~Li},
  \bibinfo{author}{C.~Lubich}, \bibinfo{title}{A convergent evolving finite
  element algorithm for {W}illmore flow of closed surfaces},
  \bibinfo{journal}{Numer. Math.} \bibinfo{volume}{149} (\bibinfo{year}{2021})
  \bibinfo{pages}{595--643}.
\bibitem[{Lee and Choi(2023)}]{LeeC23}
\bibinfo{author}{S.~Lee}, \bibinfo{author}{Y.~Choi},
  \bibinfo{title}{Curvature-based interface restoration algorithm using
  phase-field equations}, \bibinfo{journal}{PLOS ONE} \bibinfo{volume}{18}
  (\bibinfo{year}{2023}) \bibinfo{pages}{1--15}.
\bibitem[{Olischl{\"a}ger and Rumpf(2009)}]{OlischlagerR09}
\bibinfo{author}{N.~Olischl{\"a}ger}, \bibinfo{author}{M.~Rumpf},
  \bibinfo{title}{Two step time discretization of {W}illmore flow}, in:
  \bibinfo{editor}{E.R. Hancock}, \bibinfo{editor}{R.R. Martin},
  \bibinfo{editor}{M.A. Sabin} (Eds.), \bibinfo{booktitle}{Mathematics of
  Surfaces XIII}, \bibinfo{publisher}{Springer}, \bibinfo{address}{Berlin},
  \bibinfo{year}{2009}, pp. \bibinfo{pages}{278--292}.
\bibitem[{Rusu(2005)}]{Rusu05}
\bibinfo{author}{R.E. Rusu}, \bibinfo{title}{An algorithm for the elastic flow
  of surfaces}, \bibinfo{journal}{Interfaces Free Bound.} \bibinfo{volume}{7}
  (\bibinfo{year}{2005}) \bibinfo{pages}{229--239}.
\bibitem[{Seifert(1997)}]{Seifert97}
\bibinfo{author}{U.~Seifert}, \bibinfo{title}{Configurations of fluid membranes
  and vesicles}, \bibinfo{journal}{Adv. Phys.} \bibinfo{volume}{46}
  (\bibinfo{year}{1997}) \bibinfo{pages}{13--137}.
\bibitem[{Welch and Witkin(1994)}]{Welch94free}
\bibinfo{author}{W.~Welch}, \bibinfo{author}{A.~Witkin},
  \bibinfo{title}{Free-form shape design using triangulated surfaces}, in:
  \bibinfo{booktitle}{Proceedings of the 21st annual conference on Computer
  graphics and interactive techniques}, pp. \bibinfo{pages}{247--256}.
\bibitem[{Willmore(1993)}]{Willmore93}
\bibinfo{author}{T.J. Willmore}, \bibinfo{title}{Riemannian geometry},
  \bibinfo{publisher}{Oxford University Press}, \bibinfo{year}{1993}.

\end{thebibliography}
\end{document}